\documentclass[a4paper,12pt]{article}

\usepackage[left=2.50cm, right=2.50cm, top=2.50cm, bottom=2.50cm]{geometry}
\usepackage{amsmath,amsthm,amssymb,amsfonts}
\usepackage{graphicx}
\usepackage{subfigure}
\usepackage{bm}
\usepackage{multirow}
\usepackage{booktabs}
\usepackage{longtable}
\usepackage{dsfont}
\usepackage{array}
\usepackage{indentfirst}
\usepackage{enumitem}
\setlist[itemize]{noitemsep}
\usepackage{abstract}

\usepackage{titlesec}
\usepackage[ruled]{algorithm2e}
\usepackage{algpseudocode}

\usepackage{titling}
\usepackage[round]{natbib}
\usepackage{hyperref}
\hypersetup{colorlinks, unicode}
\usepackage{tkz-euclide,pgfplots}
\usetikzlibrary{positioning}
\pgfplotsset{compat=1.17}

\AtBeginDocument{\sloppy\hfuzz=20pt}

\newtheorem{lemma}{Lemma}
\newtheorem{theorem}{Theorem}

\newtheorem{assumption}{Assumption}
\renewcommand{\theassumption}{\Alph{assumption}}
\newtheorem{corollary}{Corollary}
\newtheorem{definition}{Definition}
\newtheorem{remark}{Remark}

\newlist{subasm}{enumerate}{1}
\setlist[subasm]{
    label=(\theassumption\arabic*),
    ref=\theassumption\arabic*,
    itemsep=0pt, topsep=0pt,
    partopsep=0pt, parsep=0pt,
}

\usepackage{fancyhdr}
\title{Byzantine-Robust Distributed Sparse Softmax Regression}

\author{Yuxuan Wang\textsuperscript{a}\thanks{E-mail: \texttt{wangyuxuan@zju.edu.cn}}\quad  Yu Sun\textsuperscript{b}\thanks{ E-mail: \texttt{sunyu5968@163.com}} \quad Bowen Li\textsuperscript{c}\thanks{E-mail: \texttt{bornwhu@gmail.com}} \quad Kangqiang Li\textsuperscript{c}\thanks{Corresponding author. E-mail: \texttt{11935023@zju.edu.cn}} \\[1ex]
\small \textsuperscript{a}School of Mathematical Sciences, Zhejiang University, Hangzhou, Zhejiang 310058, China \\
\small \textsuperscript{b}Economic Information Center of China National Tobacco Corporation, Beijing 100045, China \\
\small \textsuperscript{c}Information Center, Hubei Provincial Tobacco Monopoly Administration, \\
\small \phantom{\textsuperscript{c}}Wuhan, Hubei 430030, China}

\date{\today}

\begin{document}

\maketitle

\begin{abstract}
We study Byzantine-robust distributed sparse softmax regression in high-dimensional settings. We propose a communication-efficient distributed algorithm that employs robust aggregation methods---coordinate-wise trimmed mean and coordinate-wise median---to defend against Byzantine attacks. Under sub-Gaussian design assumptions with sparsity constraints, we establish that the distributed estimator achieves the  convergence rate of $O_{P}(K^{\frac{5}{2}}\sqrt{s\log(Kd)/N})$ when the Byzantine fraction $\alpha \asymp 1/\sqrt{m}$, matching the centralized oracle rate. Numerical experiments on synthetic data demonstrate the superior performance of the proposed method.
\end{abstract}


\section{Introduction}

In high-dimensional multiclass problems, only a small subset of predictors may matter for each class. Sparse softmax regression represents the conditional class probabilities through a matrix of coefficients and uses this structure for estimation and variable selection \citep{tian}. When observations are distributed across multiple workers, moving all samples to a master machine incurs communication costs. A gradient-based surrogate likelihood update instead broadcasts the current parameters, collects gradients computed from local samples, and updates the model without transferring the raw observations \citep{Communication3}.

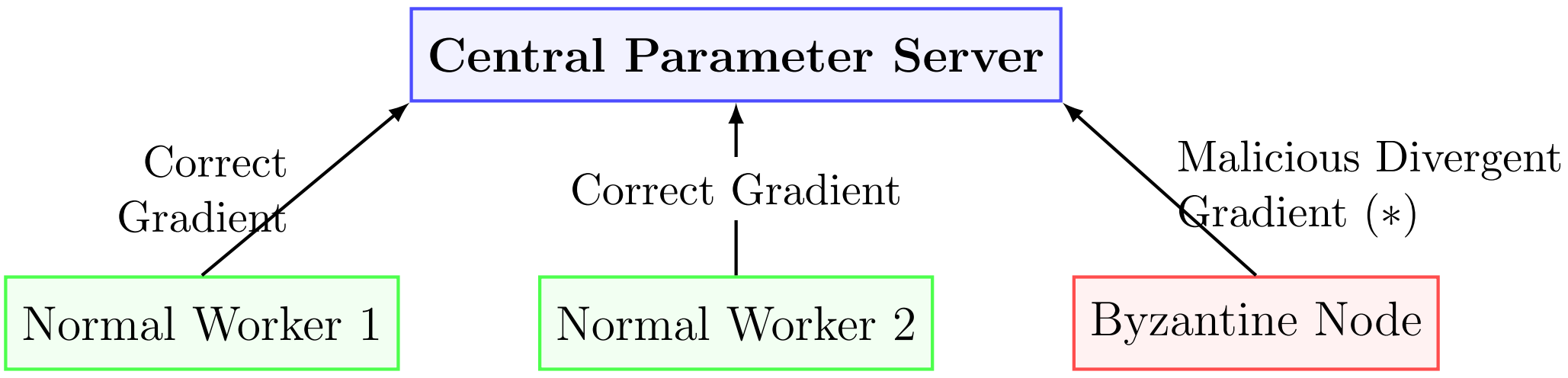
\begin{figure}[htbp]
\centering
\begin{tikzpicture}[
    node distance=1.5cm and 1.2cm,
    server/.style={rectangle, draw=blue!70, fill=blue!5, thick, minimum width=4.5cm, minimum height=0.8cm, align=center},
    worker/.style={rectangle, draw=green!70, fill=green!5, thick, minimum width=2.5cm, minimum height=0.8cm, align=center},
    byz/.style={rectangle, draw=red!70, fill=red!5, thick, minimum width=2.5cm, minimum height=0.8cm, align=center},
    arrow/.style={->, >=latex, thick}
]
    \node[server] (server) {\textbf{Central Parameter Server}};
    \node[worker] (w2) [below=of server] {Normal Worker $2$};
    \node[worker] (w1) [left=of w2] {Normal Worker $1$};
    \node[byz] (byz) [right=of w2] {Byzantine Node};

    \draw[arrow] (w1.north) -- node[left, font=\small, align=right] {Correct\\Gradient} (server.south west);
    \draw[arrow] (w2.north) -- node[fill=white, font=\small] {Correct Gradient} (server.south);
    \draw[arrow] (byz.north) -- node[right, font=\small, align=left] {Malicious Divergent\\Gradient ($\ast$)} (server.south east);
\end{tikzpicture}
\caption{Illustration of the failure mode under a Byzantine attack.}
\label{fig:byzantine_attack}
\end{figure}

This update assumes that workers report gradients computed from their local data. A Byzantine worker affected by a system failure or an attack may instead send an arbitrary gradient. Because a single arbitrarily large reported message can dominate the arithmetic mean, it can distort the direction of a model update \citep{stat,Divide}. For sparse softmax regression, the resulting error can affect both the estimated class probabilities and the selection of relevant coefficients. Figure~\ref{fig:byzantine_attack} illustrates how an arbitrary gradient can enter mean aggregation. Distributed estimation must therefore limit the effect of corrupted gradients while retaining the statistical benefit of sparsity.

For high-dimensional multinomial regression outside the distributed setting, \citet{tian} establish an error bound for the $\ell_1$-penalized estimator whose dependence on dimension is logarithmic. This provides a statistical benchmark for sparse softmax estimation. With data distributed across machines, the Communication-Efficient Surrogate Likelihood (CSL) framework of \citet{Communication3} uses local losses and exchanged gradients to approximate an update based on all observations. Related distributed sparse learning uses exchanged gradients and a shifted $\ell_1$-penalized local loss to recover sparsity-dependent accuracy under reliable communication \citep{Efficient}. CEASE \citep{Communication-efficient} develops this line of communication-efficient iterative estimation through a regularized surrogate objective. \citet{zhou2023} combine surrogate updates with robust gradient aggregation and allow convex penalties, including an $\ell_1$ penalty. Their general error bounds, however, are expressed in the ambient dimension and do not establish an $s\log(Kd)$ guarantee for sparse softmax regression.

Byzantine-robust sparse learning has been studied for robust regression, sparse SVM \citep{project1} and distributed sparse $M$-estimation \citep{Divide}. For multiclass classification, \citet{Bao} establish sparsity-dependent estimation and feature-support recovery for Byzantine-robust linear discriminant analysis under a shared-covariance Gaussian model. Their guarantees concern a group-penalized discriminant estimator and do not directly establish error or entrywise support guarantees for the conditional softmax estimator. The missing guarantee is specific: a Byzantine-robust distributed softmax estimator with a sparsity-dependent error bound and support-recovery conditions.

The proposed update substitutes a robust aggregate for the mean worker gradient in CSL. Proofs of the sparsity-dependent error rate then require bounds on the resulting gradient-correction error under Byzantine contamination. Moreover, the shared normalization in the softmax likelihood produces cross-class blocks in the Hessian. Our assumptions and support-recovery analysis therefore concern the Hessian of the full coefficient matrix.

We build on CSL and the iterative surrogate-likelihood approach of CEASE to study Byzantine-robust distributed sparse softmax regression. Our contributions are threefold. \emph{First}, we adapt the $\ell_1$-penalized surrogate update to softmax regression and use either a coordinate-wise trimmed mean or a coordinate-wise median \citep{stat} to aggregate worker gradients. \emph{Second}, we derive a sparsity-dependent Frobenius error bound with an explicit term for Byzantine contamination. When $\alpha\asymp m^{-1/2}$, the resulting rate matches the non-distributed upper bound of \citet{tian} in its dependence on $N$, $s$, and $d$. \emph{Third}, we give conditions that rule out false support coordinates and, with an additional minimum-signal assumption, guarantee exact support recovery and sign consistency.

Numerical experiments on synthetic and real data examine estimation error and predictive performance under Byzantine attacks. We also extend the error analysis to binary robit regression. Section~\ref{SecSetup} describes the distributed setting and Byzantine failures. Section~\ref{SecMain} introduces the sparse softmax estimator and states its error and support-recovery guarantees, together with the robit extension. Section~\ref{Secexper} reports the numerical experiments, and Section~\ref{SecDiscussion} discusses limitations and future directions. Proofs are given in the Appendix.

\subsection*{Notations}\label{SecDef}
For any positive integer $n$, we denote the set $\{1,2,\ldots,n\}$ by $[n]$. For two vectors $\boldsymbol{x},\boldsymbol{y} \in \mathbb{R}^{n}$, we write $\langle \boldsymbol{x},\boldsymbol{y}\rangle=\sum_{i=1}^{n}x_{i}y_{i}$ to be the inner product and for a vector $\boldsymbol{v} \in \mathbb{R}^{d}$, the support of $\boldsymbol{v}$ is defined as $\operatorname{supp}(\boldsymbol{v})=\{j: v_{j}\neq 0\}$. For $q \in [1,\infty)$, the ${\ell}_{q}$-norm is defined as $\|\boldsymbol{v}\|_{q}=\left(\sum_{j=1}^{d}|v_{j}|^{q}\right)^{1/q}$ and $\|\boldsymbol{v}\|_{\infty}=\max_{j \in [d]}|v_{j}|$. Given a subset $S$ of index set $[d]$, the projection $\boldsymbol{v}_{S}$ is defined by $\left(\boldsymbol{v}_{S}\right)_{j}=0, \operatorname{if} j \notin S $ and $ \left(\boldsymbol{v}_{S}\right)_{j}=(\boldsymbol{v})_{j},$ if  $j \in S$,  For a matrix $\boldsymbol{A}=(a_{ij}) \in \mathbb{R}^{n_{1}\times n_{2}}$, the $\|\cdot\|_{\max}$-norm of $\boldsymbol{A}$ is defined as $\|\boldsymbol{A}\|_{\max}=\max_{i \in [n_{1}],j \in [n_{2}]}|a_{ij}|$ and define $\|\boldsymbol{A}\|_{F}$  to be its Frobenius norm. For a matrix $\boldsymbol{M}$, write $|||\boldsymbol{M}|||_{\infty}:=\max_{i}\sum_{j}|M_{ij}|$ for the row-sum $\ell_{\infty}$ norm, which is to be distinguished from the entrywise norm $\|\boldsymbol{M}\|_{\max}$. 
$\boldsymbol M_{AB}$ denotes the submatrix with rows indexed by set $A$ and columns indexed by set $B$. Given two sequences $\{a_{n}\}_{n=1}^{\infty}$ and $\{b_{n}\}_{n=1}^{\infty}$, we define $a_{n}=O(b_{n})$ or $a_{n} \lesssim b_{n}$ to denote that for some positive constant $C$, $a_{n}\leq Cb_{n}$ holds for all $n$ large enough and we use the notation $a_{n} \asymp b_{n}$, if $a_{n}=O(b_{n})$ and $b_{n}=O(a_{n})$. Denote $x_n=O_{\mathbb{P}}\left(y_n\right)$ if for any $\epsilon>0, \exists C>0 \text {, such that } \mathbb{P}\left(\left|x_n / y_n\right|>C\right)<\epsilon \text { holds for all } n \geq 1$. For a differentiable function $f:\mathbb{R}^{d}\to \mathbb{R}$, we denote the gradient of $f$ by $\nabla f$. For any metric space $(\mathcal{F},\|\cdot\|)$, we denote the  $\epsilon$-cover number by $N(\epsilon, \mathcal{F},\|\cdot\|) $. 
\begin{definition}\label{def1}
 (Coordinate-wise median). For $m$ vectors of dimension $d$, $\left\{\boldsymbol{v}_i : i \in[m]\right\}$, the coordinate-wise median $\boldsymbol{g}:=\operatorname{median}\left\{\boldsymbol{v}_i: i \in[m]\right\}$ is a $d$-dimensional vector whose $k$-th element is $g_k=\operatorname{median}\left\{x_i^{(k)}: i \in[m]\right\}$, where $k \in[d]$.
\end{definition}
\begin{definition}\label{def2}
(Coordinate-wise trimmed mean). For $\beta \in\left[0, \frac{1}{2}\right)$ and $m$ $d$ dimensional vectors $\left\{\boldsymbol{v}_i : i \in[m]\right\}$, the coordinate-wise $\beta$-trimmed mean $\mathbf{g}:=\operatorname{trmean}_\beta\left\{\mathbf{x}^i: i \in[m]\right\}$ is a vector with its $k$-th coordinate being $g_k=\frac{1}{(1-2 \beta) m} \sum_{x \in U_k} x$ for each $k \in[d]$. Here $U_k$ is a subset of $\left\{x_k^1, \ldots, x_k^m\right\}$ obtained by removing the largest and smallest $\beta$ fraction of its elements.
\end{definition}
\section{Problem Setup}\label{SecSetup}
The traditional approach for recovering the sparse regression parameter $\boldsymbol{\Theta}^{*}\in \mathbb{R}^d$ in high-dimensional settings is to minimize the following optimization problem:
\begin{equation}\label{eq1}
	\boldsymbol{\Theta}^{\star}:=\underset{\boldsymbol{\Theta} \in \mathbb{R}^{d}}{\operatorname{argmin}} \{\widehat{\ell}_{N}(\boldsymbol{\Theta}; \{(y_{i},X_{i})\}_{i\in [N]})+\lambda\|\boldsymbol{\Theta}\|_{1}\}
\end{equation}
where $\lambda\|\boldsymbol{\Theta}\|_{1}$ is the penalty function to promote sparsity of the solution, and $\widehat{\ell}_{N}(\boldsymbol{\theta};\{(y_{i},X_{i})\}_{i\in [N]})$ is the empirical loss function based on $N$ data samples $\{(y_{i},X_{i})\}_{i\in [N]}$.
When both the sample size $N$ and the dimension $d$ are large, solving~\eqref{eq1} on a single machine is computationally impractical. Moreover, in practice, data are often distributed across multiple servers. Specifically, $N$ samples $\{(y_{i}^{(k)},X_{i}^{(k)})\}_{i\in [n]}^{k \in [m]}$ are evenly distributed among $m$ servers (i.e., each server contains $n=N/m$ sub-samples). 

Communication between servers incurs costs, making it economically infeasible to collect all data on a single server to solve~\eqref{eq1}. Therefore, designing a communication- and computation-efficient distributed algorithm is necessary. Based on the above considerations, in this paper we study the design of distributed regression algorithms that are resilient to Byzantine attacks, aiming to improve the statistical performance of distributed learning. We assume that a fraction $\alpha< \frac{1}{2}$ of the $m$ servers are Byzantine, while the remaining $1-\alpha$ fraction are normal. Information transmitted by Byzantine servers is denoted by $*$. To defend against Byzantine attacks, we employ the median during information aggregation.
%
\paragraph{Subsets induced by trimmed mean.} For each coordinate $(j,k)$, let $\mathcal{T} \subset [m]$ denote the set of server indices whose value at that coordinate is among the $\beta m$ largest or $\beta m$ smallest (i.e., the trimmed-out set), and let $\mathcal{U} := [m]\setminus \mathcal{T}$ denote the set of indices retained by $\operatorname{trmean}_\beta$. By construction $|\mathcal{U}|=(1-2\beta)m$.

\section{Main Results} \label{SecMain}
\subsection{Softmax Regression}
We study distributed softmax regression under Byzantine robustness. Suppose there are $N$ $i.i.d.$ samples $\{(y_{i},\boldsymbol{x}_{i})\}_{i\in [N]}$ on $m$ servers from the softmax regression model:
$$
\mathbb{P}(y=k|\boldsymbol{x})=\frac{\exp\left(\langle \boldsymbol{\theta}^{\star}_{k},\boldsymbol{x}\rangle\right)}{1+\sum\limits_{i=1}^{K-1}\exp\left(\langle \boldsymbol{\theta}^{\star}_{i},\boldsymbol{x}\rangle\right)}, \quad k \in [K-1],
$$
where $\boldsymbol{\Theta}^{\star}:=\left(\boldsymbol{\theta}_{1}^{\star},\ldots, \boldsymbol{\theta}_{K-1}^{\star}\right)\in \mathbb{R}^{d\times (K-1)}$ is the parameter matrix. This paper considers the case where $\boldsymbol{\Theta}^{\star}$ is sparse and high-dimensional, i.e., $\max\limits_{k}\|\boldsymbol{\theta}^{\star}_{k}\|_{0}\leq s\ll d$. To estimate the sparse parameter matrix $\boldsymbol{\Theta}^{\star}$, a Lasso penalty on the negative log-likelihood function is imposed and the idealized approach is the centralized estimator, i.e. minimizing the following optimization problem:
    $$
	\widehat{\boldsymbol{\Theta}} :=\underset{\boldsymbol{\Theta} \in \mathbb{R}^{d\times (K-1)}}{\operatorname{argmin}}\left\{ \widehat{L}(\boldsymbol{\Theta})+\lambda\|\boldsymbol{\Theta}\|_{1,1} \right\},
$$
where
$\widehat{L}(\boldsymbol{\Theta}):=\frac{1}{N} \sum_{i=1}^N \log \left(1+\sum_{k=1}^{K-1} \exp \left(\boldsymbol{\theta}_{k}^{\top} \boldsymbol{x}_i\right)\right)-\frac{1}{N} \sum_{k=1}^{K-1}\sum_{i=1}^N  y_i^{(k)}\boldsymbol{\theta}_{k}^{\top} \boldsymbol{x}_i$ and $y_{i}^{(k)}=\mathbf{1}(y_{i}=k).$
\citet{tian} shows that the centralized approach achieves the statistical error rate of $O_{\mathbb{P}} \left(K^{\frac{5}{2}}\sqrt{s \log (Kd)/N}\right)$ under Frobenius norm, however, it is impractical due to expensive communication.\par 
To solve this problem, \citet{Communication3} proposed the following  empirical loss function with reduced communication cost:
\begin{equation}\label{eq:csl-surrogate}
\widetilde{L}(\boldsymbol{\beta})=\widehat{L}_1(\boldsymbol{\beta})-\boldsymbol{\beta}^{\top}\left(\nabla \widehat{L}_1(\tilde{\boldsymbol{\beta}})-\frac{1}{m}\sum_{i\in [m]}\nabla \widehat{L}_{i}(\tilde{\boldsymbol{\beta}})\right)
\end{equation}
where $\tilde{\boldsymbol{\beta}}$ is an initialize estimator. When Byzantine failures exist, the mean-gradient correction in~\eqref{eq:csl-surrogate} is sensitive to incorrect gradient values. The gradients $\{\nabla \widehat{L}_{i}(\tilde{\boldsymbol{\beta}}): i = 1,\ldots,m\}$ received by the central machine are not always reliable. Therefore, based on coordinate-wise median and truncated mean, we propose the following Byzantine-robust distributed algorithm. The communication cost of Algorithm~\ref{Al2} is $O(dmK)$, since each server only transmits $\widehat{\boldsymbol{\Theta}}_{0}$ and $\widehat{\boldsymbol{G}}^{(\ell)}\left(\widehat{\boldsymbol{\Theta}}_{0}\right)$, and $K\ll n$. This is much smaller than the communication cost $O(dmn)$ of traditional methods.

\begin{algorithm}[!ht]
	\caption{Byzantine-Robust Distributed Algorithm} 
	\label{Al2}
	\KwIn{Initialize estimator $\widehat{\boldsymbol{\Theta}}_{0}$, algorithm parameter $\beta$ (for Option I) and number of iteration $T$.}
	\For{$t=0,1,2,\ldots,T-1$}{
		\For{$\ell=1,2,3,\ldots,m$}{
			Each machine $\ell$: evaluate local gradient
			$$
			\widehat{\boldsymbol{G}}^{(\ell)}\left(\widehat{\boldsymbol{\Theta}}_{t}\right) \leftarrow \begin{cases}\nabla{\widehat{L}}_{\ell}(\widehat{\boldsymbol{\Theta}}_{t}) & \text { normal worker machines, } \\ * & \text { Byzantine machines, }\end{cases}
			$$
			and send to the 1st machine.
            }
		 The 1st machine: evaluates aggregate gradients through robust methods: $$\nabla \widetilde{\mathcal{L}}(\widehat{\boldsymbol{\Theta}}_{t})\leftarrow \begin{cases}\operatorname{trmean}_{\beta}\left\{\widehat{\boldsymbol{G}}^{(\ell)}\left(\widehat{\boldsymbol{\Theta}}_{t}\right): \ell \in [m]\right\} & \text { Option \uppercase\expandafter{\romannumeral1}, } \\ \operatorname*{median}_{\ell \in [m]}\left\{\widehat{\boldsymbol{G}}^{(\ell)}\left(\widehat{\boldsymbol{\Theta}}_{t}\right)\right\} & \text { Option \uppercase\expandafter{\romannumeral2},}\end{cases}
		$$
		solves the following regularized optimization problem:
		\begin{equation}\label{eq3}
		\widehat{\boldsymbol{\Theta}}_{t+1}:=\underset{\boldsymbol{\Theta} \in \mathbb{R}^{d\times (K-1)}}{\operatorname{argmin}} \left\{\widehat{L}_{1}(\boldsymbol{\Theta})-
		\left\langle \nabla \widehat{L}_{1}(\widehat{\boldsymbol{\Theta}}_{t})-\nabla \widetilde{\mathcal{L}}(\widehat{\boldsymbol{\Theta}}_{t}),\boldsymbol{\Theta}\right\rangle+\lambda_{t+1}\|\boldsymbol{\Theta}\|_{1,1}
		\right\}
	    \end{equation}
		and then broadcasts it to the other $m-1$ machines.
	}\KwOut{$\widehat{\boldsymbol{\Theta}}_{T}$.}
\end{algorithm}
\begin{assumption}\label{assumption1}
\begin{subasm}
	\item $\boldsymbol{x}$ is a $d$-dimensional sub-Gaussian random vector with $\mathbb{E} [\boldsymbol{x}]=\boldsymbol{0}_{d\times 1},\lambda_{\min}(\Sigma)\geq c>0, \Sigma:=\mathbb{E} \left[\boldsymbol{x}\boldsymbol{x}^{\top}\right]$.\label{a1}
	\item $\max\limits_{k}\left\|\boldsymbol{\theta}^{\star}_{k}\right\|_{0}\leq s\ll d$, $\max\limits_{k}\left\|\boldsymbol{\theta}^{\star}_{k}\right\|_{2} \leq C_{1}$.\label{a2}
\end{subasm}
\end{assumption}
Under condition~\ref{a1}, $\left\|\boldsymbol{x}\right\|_{\max} \leq c_d:=C \sqrt{\log d}$ with probability $1-d^{-C}$. It can be shown that these conditions are all easily satisfied.\par
Let
$$
S := \bigl\{(j,k)\in[d]\times[K-1]:\ \Theta^{\star}_{j,k}\neq 0\bigr\},\qquad |S|\leq (K-1)s,
$$
denote the true support.
The Hessian matrix of the local empirical softmax loss $\widehat{L}_1$ can write as
$$
\nabla^2\widehat L_1(\boldsymbol\Theta)=\frac1n\sum_{i=1}^n \boldsymbol A_i(\boldsymbol\Theta)\otimes \boldsymbol x_i\boldsymbol x_i^\top,
\quad
\boldsymbol A_i(\boldsymbol\Theta)=\operatorname{diag}(\boldsymbol p_i)-\boldsymbol p_i\boldsymbol p_i^\top\in\mathbb R^{(K-1)\times(K-1)},
$$
and let $\boldsymbol{H}^{\star}:=\mathbb{E}\bigl[\boldsymbol{A}(\boldsymbol{\Theta}^{\star})\otimes\boldsymbol{x}\boldsymbol{x}^{\top}\bigr]\in\mathbb{R}^{d(K-1)\times d(K-1)}$ be the population Hessian, with $\boldsymbol{p}=(p_{1}(\boldsymbol{x},\boldsymbol{\Theta}),\dots,p_{K-1}(\boldsymbol{x},\boldsymbol{\Theta}))^{\top}$. Here $\boldsymbol{H}^{\star}$ is the Hessian with respect to the vectorized parameter $\operatorname{vec}(\boldsymbol{\Theta})\in\mathbb{R}^{d(K-1)}$. 

\begin{assumption}\label{assumption2}
There exist constants $\gamma\in(0,1]$, $C_{\min}>0$, and $D_{\infty}<\infty$ such that
\begin{subasm}
	\item $|||\boldsymbol{H}^{\star}_{S^{c}S}(\boldsymbol{H}^{\star}_{SS})^{-1}|||_{\infty}\leq 1-\gamma$;\label{a2cond1}
	\item $\lambda_{\min}\bigl(\boldsymbol{H}^{\star}_{SS}\bigr)\geq C_{\min}$;\label{a2cond2}
	\item $|||(\boldsymbol{H}^{\star}_{SS})^{-1}|||_{\infty}\leq D_{\infty}$.\label{a2cond3}
\end{subasm}
\end{assumption}
Assumption~\ref{assumption2} extends the incoherence and dependency conditions of \citet{Ravikumar2010} from the single-class Fisher information matrix to the population Hessian $\boldsymbol{H}^{\star}$ arising in multi-class softmax regression.

\begin{figure}[!ht]
	\centering
	\includegraphics[width=4.7in,height=2.7in]{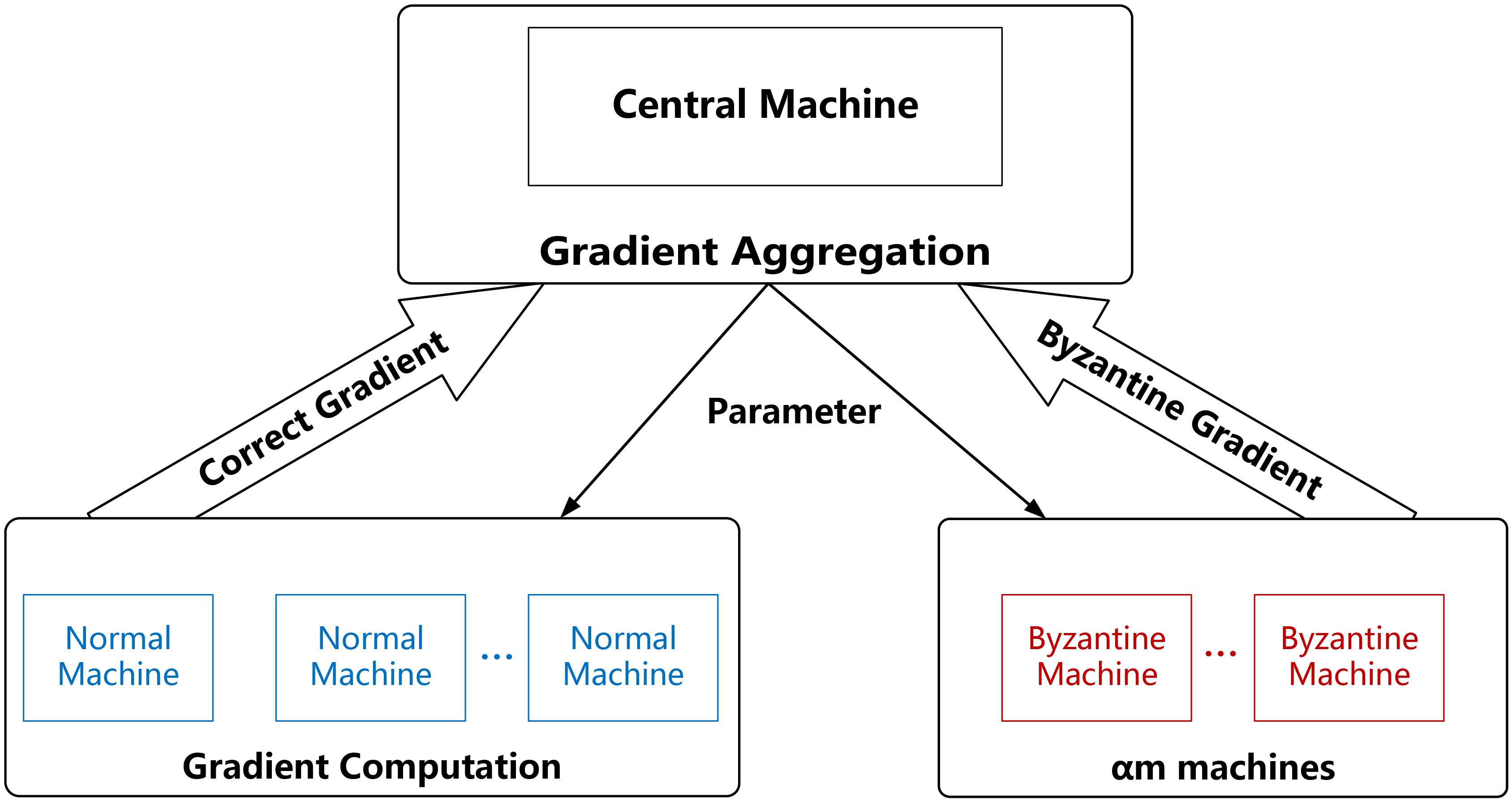}
	\caption{The illustration of Algorithm~\ref{Al2}. }
	\label{fig:alg}
\end{figure}
\begin{theorem}\label{theorem1}
	Under Assumptions~\ref{assumption1} and~\ref{assumption2}, choosing
    $$\lambda_{1}:= C_{1}\sqrt{\frac{\log (Kd)}{N}}+C_{2}\log(d)^3K^6s^2\frac{\log (Kd)}{n}+C_{3}\frac{2\beta+\alpha}{1-2\beta}\sqrt{\frac{\log(Kd)}{n}},$$
 provided that $\log(d)^6K^{12}s^4\log(Kd)m\lesssim n$, the estimator $\widehat{\boldsymbol{\Theta}}_{1}$ produced by Option $\uppercase\expandafter{\romannumeral1}$ satisfies the following bounds with probability at least $1-cmd^{-C}$,
	$$
	\begin{aligned}
	&\left\|\widehat{\boldsymbol{\Theta}}_{1}-\boldsymbol{\Theta}^{\star}\right\|_{F} \leq C K^{\frac{5}{2}}\left(\sqrt{\frac{s\log(Kd)}{N}}+\frac{2\beta\!+\!\alpha}{1-2\beta}\sqrt{\frac{s\log(Kd)}{n}}\right),\\ &
	\left\|\widehat{\boldsymbol{\Theta}}_{1}-\boldsymbol{\Theta}^{\star}\right\|_{1,1} \leq C\,s\,K^{3}\left(\sqrt{\frac{\log(Kd)}{N}}+\frac{2\beta\!+\!\alpha}{1-2\beta}\sqrt{\frac{\log(Kd)}{n}}\right).
	\end{aligned}
	$$
    Furthermore, by math induction, we also have that for $\forall t\geq 1$,
    $$\left\|\widehat{\boldsymbol{\Theta}}_{t+1}-\boldsymbol{\Theta}^{\star}\right\|_{F} \lesssim_{\mathbb{P}} K^{\frac{5}{2}}\left(\sqrt{\frac{s\log(Kd)}{N}}+\frac{2\beta+\alpha}{1-2\beta}\sqrt{\frac{s\log(Kd)}{n}}\right),$$
    when $\lambda_{t+1} \asymp \sqrt{\frac{\log (Kd)}{N}}+\frac{2\beta+\alpha}{1-2\beta}\sqrt{\frac{\log(Kd)}{n}}$.
\end{theorem}

\begin{remark}\label{remark1}
	If $c\alpha\geq \beta\geq \alpha$, the estimator under Option $\uppercase\expandafter{\romannumeral1}$ achieves the convergence rate
	$$
	\left\|\widehat{\boldsymbol{\Theta}}_{T}-\boldsymbol{\Theta}^{\star}\right\|_{F} = O_{\mathbb{P}} \left(K^{\frac{5}{2}}\sqrt{\frac{s \log(Kd)}{N}}+\alpha K^{\frac{5}{2}}\sqrt{\frac{s \log(Kd)}{n}}\right)
	$$
	This is consistent with the result proved by \citet{stat}: when there are no Byzantine servers in the algorithm, i.e., $\alpha=0$, Theorem~\ref{theorem1} implies that the distributed estimator $\widehat{\boldsymbol{\Theta}}_{T}$ achieves the optimal convergence rate (\citet{tian}); if $\alpha \asymp \frac{1}{\sqrt{m}}$, then $\widehat{\boldsymbol{\Theta}}_{T}$ can still maintain the statistical error rate of $O_{\mathbb{P}} \left(K^{\frac{5}{2}}\sqrt{s \log (Kd)/N}\right)$.
\end{remark}
\begin{remark}\label{remark2}
	Compared to \citet{zhou2023}, the estimator produced by our algorithm has the following advantages. Theoretically, the upper bound of our estimator is superior to the $O_{\mathbb{P}}\left(\sqrt{\frac{d}{N}}+\alpha\sqrt{\frac{d}{n}}\right)$ bound in \citet{zhou2023}.
\end{remark}

\subsubsection{Sparsity recovery}
To streamline subsequent statements, define
\begin{equation}\label{eq:rhoN}
	\rho_{N} \;:=\; \sqrt{\frac{\log(Kd)}{N}} \;+\; \frac{\sqrt{1-\alpha}}{1-2\beta}\sqrt{\frac{\log(Kd)}{N}} \;+\; \frac{2\beta+\alpha}{1-2\beta}\sqrt{\frac{\log(Kd)}{n}}.
\end{equation}
Under an additional minimum-signal condition, the conclusion of Theorem~\ref{theorem2} can be strengthened to exact support recovery and sign consistency, since the support-inclusion part is obtained as the first half of the argument for Theorem~\ref{theorem2prime}.

\begin{theorem}\label{theorem2}
	Suppose Assumptions~\ref{assumption1} and~\ref{assumption2} hold, and that the regularization parameter satisfies $\lambda\geq \tfrac{8}{\gamma}\,\rho_{N}$. If $n\gtrsim \gamma^{-2}\operatorname{poly}(K)\,s^{2}\log(Kd)$, then there exist constants $c_{1},c_{2}>0$ independent of $n,m,d,K$ such that for $\forall t\geq 0$,
	$$
	\mathbb{P}\bigl(\operatorname{supp}(\widehat{\boldsymbol{\Theta}}_{t})\subseteq S\bigr)\;\geq\;1-c_{1}(Kd)^{-c_{2}}.
	$$
\end{theorem}

\subsection{The Robit Model}
As an extension of our main results, we now show that Algorithm~\ref{Al2} also applies to the binary robit regression model. The robit model replaces the logistic link in logit regression with the cumulative distribution function of a Student-$t$ distribution, which is more robust than the logit and probit links when the response is contaminated by label-flipping noise; see \citet{liu} for the original proposal and \citet{tian} for the high-dimensional treatment. Specifically,
$$
\mathbb{P}\left(y=1 \mid \boldsymbol{x}, \boldsymbol{\theta}\right)=1-\mathbb{P}\left(y=0 \mid \boldsymbol{x}, \boldsymbol{\theta}\right)=F_\nu\left(\boldsymbol{x}^{\top} \boldsymbol{\theta}\right),
$$
where $F_\nu(x)$ denotes the cdf of a Student-$t$ random variable with mean zero, unit scale, and $\nu$ degrees of freedom, with density
$$
f_\nu(x) = \frac{\Gamma((\nu+1)/2)}{(\pi \nu)^{1/2}\,\Gamma(\nu/2)\,\bigl(1+x^{2}/\nu\bigr)^{(\nu+1)/2}}, \qquad x\in(-\infty,\infty).
$$
As $\nu\to\infty$, the robit model reduces to the probit regression model.

\begin{corollary}\label{corollary1}
	Suppose that the sub-Gaussian design and sparsity conditions in Assumption~\ref{assumption1} hold, with the softmax local strong convexity condition replaced by the analogous condition for the robit log-likelihood in a neighborhood of $\boldsymbol{\theta}^{\star}$. Replace the local gradient in Algorithm~\ref{Al2} by the gradient of the robit log-likelihood. Then, provided that $\log(d)^{6} s^{4} \log(d)\,m \lesssim n$, the output estimator $\widehat{\boldsymbol{\theta}}_{1}$ satisfies, with high probability,
	$$
	\bigl\|\widehat{\boldsymbol{\theta}}_{1}-\boldsymbol{\theta}^{\star}\bigr\|_{2} \;\leq\; C\!\left(\sqrt{\frac{s\log d}{N}}+\frac{2\beta+\alpha}{1-2\beta}\sqrt{\frac{s\log d}{n}}\right),
	$$
	where $C$ is a constant independent of $n,m,d$ and depends on the degrees of freedom $\nu$ and the smallest eigenvalue of the design covariance $\boldsymbol{\Sigma}$.
\end{corollary}

\section{Numerical Experiments} \label{Secexper}
In this section, we conduct extensive numerical experiments on synthetic data. We use the log estimation error $\log\left\|\widehat{\boldsymbol{\Theta}}_{t}-\boldsymbol{\Theta}^{\star}\right\|_{F}$ to demonstrate the superior performance of our algorithm.
$$
\widehat{\boldsymbol{\theta}}_{k}^{(j+1)}=\widehat{\boldsymbol{\theta}}_{k}^{(j)}-\left(\frac{1}{n} \mathbf{X}^{\top} \mathbf{W}_k^{(j)} \mathbf{X}\right)^{-1}\left[\mathbf{X}^{\top}\left(\mathbf{P}^{(j)}-\mathbf{Y}\right)/n-\nabla \widehat{L}_{1}\left(\widehat{\boldsymbol{\Theta}}_0\right)+\nabla \widetilde{\mathcal{L}}\left(\widehat{\boldsymbol{\Theta}}_0\right)\right]_{:,k},
$$
where
\begin{align*}
\bm{\Phi}&:=\mathbf{X}^{\top}\left(\mathbf{P}^{(j)}-\mathbf{Y}\right)/n-\nabla \widehat{L}_{1}\left(\widehat{\boldsymbol{\Theta}}_0\right)+\nabla \widetilde{\mathcal{L}}\left(\widehat{\boldsymbol{\Theta}}_0\right),\\
\mathbf{P}^{(j)}&:=e^{\left(\boldsymbol{\theta}_{k}^{(j)}\right)^{\top} \boldsymbol{x}_i} /\left(1+\textstyle\sum_{k=1}^{K}e^{\left(\boldsymbol{\theta}_{k}^{(j)}\right)^{\top} \boldsymbol{x}_i}\right),\quad (\mathbf{Y})_{i,k}=\mathbf{1}(y_{i}=k),\\
\mathbf{W}_k^{(j)} &\in \mathbb{R}^{n\times n}:=e^{\left(\boldsymbol{\theta}_{k}^{(j)}\right)^{\top} \boldsymbol{x}_i} /\left(1+\textstyle\sum_{k=1}^{K}e^{\left(\boldsymbol{\theta}_{k}^{(j)}\right)^{\top} \boldsymbol{x}_i}\right)^2.
\end{align*}

\subsection{Effect of the number of servers \texorpdfstring{$m$}{m} (no Byzantine workers)}\label{sec:exp_m}
We first consider distributed sparse softmax regression in the Byzantine-free regime to isolate the effect of the number of servers $m$ and the local sample size $n=N/m$ on the estimation error. We set $K=4$, $s=8$, $d\in\{40,60,80,100\}$, total sample size $N=44100$, and regularization parameter $\lambda=4\sqrt{\log(d)/N}$. For each $k$, the support $S_{k}$ of $\boldsymbol{\theta}_{k}^{\star}$ is generated by sampling $s$ indices without replacement from $[d]$, and $[\boldsymbol{\theta}_{k}^{\star}]_{j}=1$ for $j\in S_{k}$. All observations $\{(y_{i},\boldsymbol{x}_{i})\}_{i=1}^{N}$ are i.i.d.\ with $\boldsymbol{x}_{i}\sim N(\boldsymbol{0}_{d\times 1},I_{d})$. Reported estimation errors are averaged over 30 independent replications.\par
\begin{figure}[h]
	\centering
	\subfigure{
		\begin{tikzpicture}
			\pgfplotsset{every axis legend/.append style={at={(0.5,0.5)},anchor=south},every axis y label/.append style={at={(-0.08,0.5)}}}
			\begin{axis}[xlabel = $\log(m)$,ylabel = {$\log(error)$}, ymajorgrids=true,legend columns=4,legend pos=north west,title style={font=\small},legend style={font=\small},font=\small,height=6.5cm, width=14cm]
				\addplot[
				color=purple,
				mark=+]
				coordinates {
					(0,-1.617656)(0.69314718,-1.5999926)(1.098612,-1.60162)(1.609438,-1.585932)(2.302585,-1.59379)(2.639057,-1.573304)(2.99573,-1.558335)(3.555348,-1.520963)(3.912023,-1.510289)(4.1431347,-1.47498)};
				\addplot[
				color=lime,
				mark=triangle]
				coordinates {
                    (0,-1.705585)(0.69314718,-1.689932388)(1.098612,-1.68549)(1.609438,-1.690637)(2.30258,-1.6772)(2.639057,-1.676492)(2.99573,-1.664612)(3.555348,-1.6486487)(3.912023,-1.620297)(4.1431347,-1.598578)};
				\addplot[
				color=cyan,
				mark=o]
				coordinates {
					(0,-1.7836365)(0.69314718,-1.777773)(1.098612,-1.7854386)(1.609438,-1.79720)(2.302585,-1.798623)(2.639057,-1.797023)(2.99573,-1.782673)(3.555348,-1.768995)(3.912023,-1.750059)(4.1431347,-1.739083)};
				\addplot[
				color=red,
				mark=square]
				coordinates {
					(0,-1.923039)(0.69314718,-1.9041647)(1.098612,-1.920745)(1.609438,-1.93139)(2.302585,-1.934812)(2.639057,-1.936464)(2.99573,-1.927127)(3.555348,-1.93523)(3.912023,-1.923075)(4.1431347,-1.8951725)};
				\legend{$d=100$,$d=80$, $d=60$,$d=40$}
			\end{axis}
	\end{tikzpicture}}
	\caption{Trend of estimation error with varying $m$.}\label{fig1}
\end{figure}
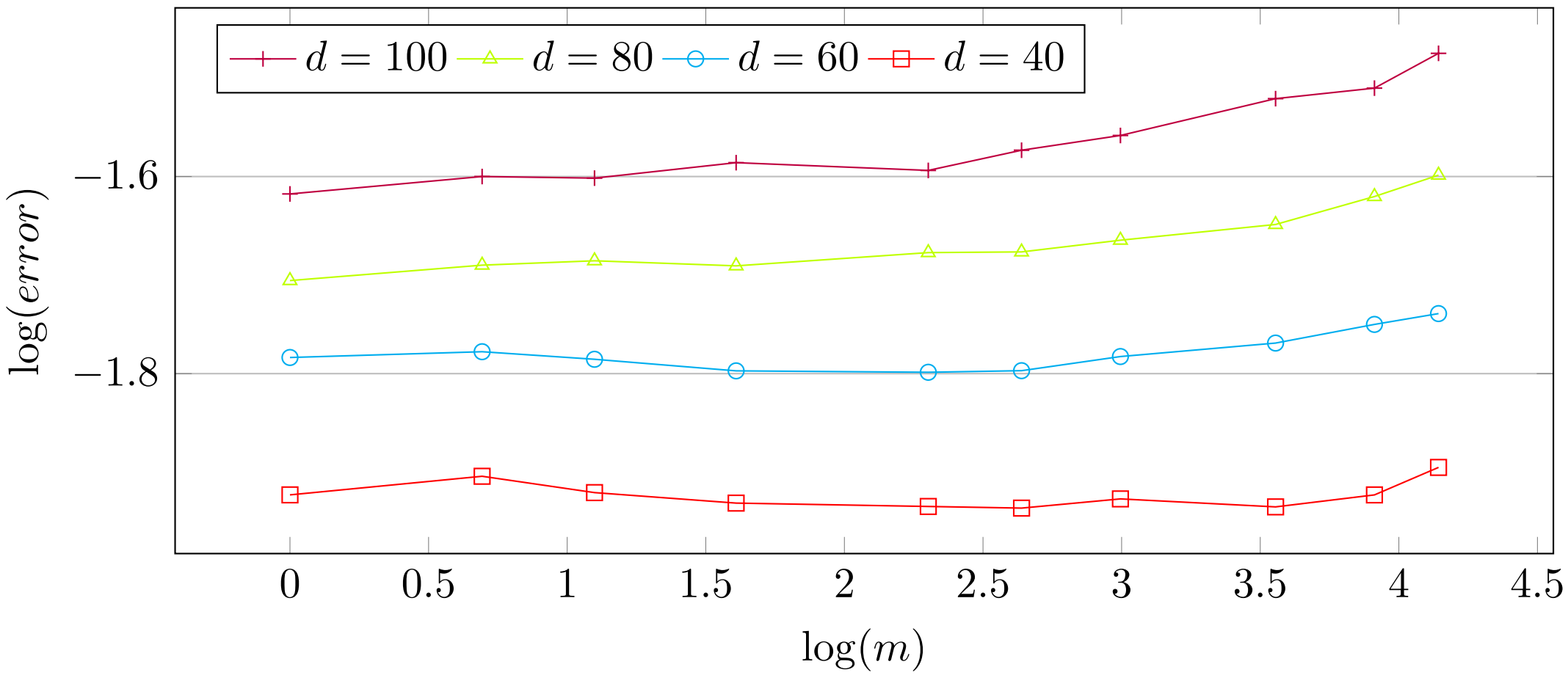

Figure~\ref{fig1} displays how the estimation error varies with $m$ for different choices of $d$. As $m$ increases (so that the local sample size $n$ decreases), the error grows slowly, and remains lower for smaller $d$. The mild growth confirms the theoretical prediction that the optimal rate is preserved when $\alpha=0$.

\subsection{Convergence under Byzantine attacks}\label{sec:exp_convergence}
\begin{figure}
	\centering
	\subfigure{
		\begin{tikzpicture}
			\tikzstyle{every node}=[font=\tiny]
			\pgfplotsset{every axis legend/.append style={
					at={(0.5,1.03)},
					anchor=south},
				every axis y label/.append style={at={(-0.08,0.5)}}}
			\begin{axis}[title={large $n=3000$},  
				xtick =data, xlabel = $T$,ylabel = {$\log(\text{error})$}, legend columns=1,legend pos=north east,legend style={font=\footnotesize},font=\small,height=5.5cm, width=5.8cm]
				\addplot[
				color=blue,
				mark=triangle
				]
				coordinates {
					(0,-0.2591628)(1, -1.6372718)(2,-1.7905198)(3,-1.8048440)(4, -1.8046880)(5, -1.8051689)(6,-1.8050642
					)(7, -1.8051074) };
				\addplot[
				color=red,
				mark=square
				]
				coordinates {
					(0,-0.260140)(1, -1.622002)(2, -1.738982)(3, -1.744048 )(4,-1.743857)(5,  -1.743837) (6,  -1.743817)(7,  -1.743819)  };
				\addplot[
				color=green,
				mark=o
				]
				coordinates {
				(0,-1.905199)(1,-1.905199)(2,-1.905199)(3,-1.905199)(4,-1.905199)(5, -1.905199) (6, -1.905199)(7, -1.905199) 
				};
				
				\legend{CSL-tri, CSL-med, Centralize}
			\end{axis}
	\end{tikzpicture}}
	\subfigure{
		\begin{tikzpicture}
			\tikzstyle{every node}=[font=\tiny]
			\pgfplotsset{every axis legend/.append style={
					at={(0.5,1.03)},
					anchor=south},
				every axis y label/.append style={at={(-0.08,0.5)}}}
			\begin{axis}[title={moderate $n=1500$},
				xtick =data,xlabel = $T$,ylabel = {$\log(\text{error})$},legend columns=1,legend pos=north east,legend style={font=\footnotesize},font=\small,height=5.5cm, width=5.8cm]
				\addplot[
				color=blue,
				mark=triangle]
				coordinates
				{(0,0.1418937)(1, -1.2110663)(2,-1.6817264)(3, -1.7767255 )(4,-1.8017990)(5, -1.8051599)(6, -1.8077677)(7, -1.8077876)(8,-1.8081845)}; 
				\addplot[
				color=red,
				mark=square
				]
				coordinates {
					(0,0.1506904)(1, -1.1601309)(2,  -1.5716026)(3,  -1.6268048 )(4,-1.6432671)(5,  -1.6437677)(6,-1.6448993)(7, -1.6449832 )(8,-1.6449139)};
				\addplot[
				color=green,
				mark=o
				]
				coordinates {
					(0,-1.905199)(1,-1.905199)(2,-1.905199)(3,-1.905199)(4,-1.905199)(5, -1.905199) (6, -1.905199)(7, -1.905199)(8, -1.905199) 
				};
				\legend{CSL-tri, CSL-med, Centralize}
			\end{axis}
	\end{tikzpicture}}
	\subfigure{
		\begin{tikzpicture}
			\tikzstyle{every node}=[font=\tiny]
			\pgfplotsset{every axis legend/.append style={
					at={(0.5,1.03)},
					anchor=south},
				every axis y label/.append style={at={(-0.08,0.5)}}}
			\begin{axis}[title={small $n=1000$},
				xtick =data,xlabel = $T$,ylabel = {$\log(\text{error})$},legend columns=1,legend pos=north east,legend style={font=\footnotesize},font=\small,height=5.5cm, width=5.8cm]
				\addplot[
				color=blue,
				mark=triangle
				]
				coordinates {
					(0,0.4725553)(2,-1.1992390)(4,-1.6032166)(6,-1.7228125)(8, -1.7606501)(10, -1.7737700)(12,-1.7791335)(14,-1.7818624)
				};
				\addplot[
				color=red,
				mark=square
				]
				coordinates {(0,0.4178577)(2,-1.2551235)(4,-1.5464211)(6,-1.5825709)(8, -1.5883024)(10, -1.5893212)(12,-1.5892190)(14,-1.5889740)};
			\addplot[
			color=green,
			mark=o
			]
			coordinates {
				(0,-1.905199)(2, -1.905199)(4, -1.905199)(6, -1.905199)(8,-1.905199)(10, -1.905199) (12, -1.905199)(14,-1.905199)
			};
				\legend{CSL-tri, CSL-med, Centralize}
			\end{axis}
	\end{tikzpicture}}
	\subfigure{
		\begin{tikzpicture}
			\tikzstyle{every node}=[font=\tiny]
			\pgfplotsset{every axis legend/.append style={
					at={(0.5,1.03)},
					anchor=south},
				every axis y label/.append style={at={(-0.08,0.5)}}}
			\begin{axis}[title={small $n$, $t_{2.1}$ noises},
				xtick =data,xlabel = $T$,ylabel = {$\log(\text{error})$},legend columns=1,legend pos=north east,legend style={font=\footnotesize},font=\small,height=5.5cm, width=5.8cm]
				\addplot[
				color=blue,
				mark=triangle]
				coordinates {(0,-0.275428)(1, -1.655247)(2, -1.713123)(3, -1.762207)(4, -1.754772)(5, -1.739105)(6, -1.745622 )(7,-1.742333)
				};
				\addplot[
				color=red,
				mark=square
				]
				coordinates {
					(0, -0.2603005)(1, -1.5258662)(2, -1.6231587)(3, -1.6136673 )(4,-1.6317318)(5,  -1.6107643)(6,  -1.6233339)(7, -1.6466305 )
				};
				\addplot[
				color=green,
				mark=o
				]
				coordinates {
					(0,-1.905199)(1,-1.905199)(2,-1.905199)(3,-1.905199)(4,-1.905199)(5, -1.905199) (6, -1.905199)(7, -1.905199)
				};
				\legend{CSL-tri, CSL-med, Centralize}
			\end{axis}
	\end{tikzpicture}}
\subfigure{
	\begin{tikzpicture}
		\tikzstyle{every node}=[font=\tiny]
		\pgfplotsset{every axis legend/.append style={
				at={(0.5,1.03)},
				anchor=south},
			every axis y label/.append style={at={(-0.08,0.5)}}}
		\begin{axis}[title={small $n$, $\alpha=0.1$},
			xtick =data,xlabel = $T$,ylabel = {$\log(\text{error})$},legend columns=1,legend pos=north east,legend style={font=\footnotesize},font=\small,height=5.5cm, width=5.8cm]
			\addplot[
			color=blue,
			mark=triangle]
			coordinates {(0, 0.1415646)(1, -1.3890854)(2, -1.6790708)(3,  -1.6695586)(4, -1.6457325)(5,  -1.6443115)(6, -1.6644315)(7, -1.6582512 )
			};
			\addplot[
			color=red,
			mark=square
			]
			coordinates {
				(0,0.1345881)(1, -1.2741271)(2, -1.5541469)(3, -1.5512853 )(4,-1.5694008)(5, -1.5847879)(6, -1.5989270 )(7,-1.5952376)
			};
			\addplot[
			color=green,
			mark=o
			]
			coordinates {
					(0,-1.905199)(1,-1.905199)(2,-1.905199)(3,-1.905199)(4,-1.905199)(5, -1.905199) (6, -1.905199)(7, -1.905199)
			};
		\legend{CSL-tri, CSL-med, Centralize}
		\end{axis}
\end{tikzpicture}}
	\subfigure{
		\begin{tikzpicture}
			\tikzstyle{every node}=[font=\tiny]
			\pgfplotsset{every axis legend/.append style={
					at={(0.5,1.03)},
					anchor=south},
				every axis y label/.append style={at={(-0.08,0.5)}}}
			\begin{axis}[title={small $n$, $\alpha=0.1$},
				xtick =data,xlabel = $T$,ylabel = {$\log(\text{error})$},legend columns=1,legend pos=north east,legend style={font=\footnotesize},font=\small,height=5.5cm, width=5.8cm]
				\addplot[
				color=blue,
				mark=triangle
				]
				coordinates {
					(0,0.4263185)(1,-0.9433951)(2, -1.4475212)(3, -1.5980405 )(4,-1.6082455)(5, -1.6243227)(6,-1.6026208)(7, -1.6140412 )(8,-1.6201030)(9, -1.6239244)(10, -1.6301644 )
				};
				\addplot[
				color=red,
				mark=square
				]
				coordinates {
					(0, 0.4461263)(1, -0.8684299)(2, -1.3922438)(3, -1.5002291 )(4,-1.5419812)(5, -1.5671543
					)(6, -1.5326918)(7, -1.5572740)(8, -1.5536989)
					(9, -1.5466822)(10, -1.5463339)
				};
				\addplot[
				color=green,
				mark=o
				]
				coordinates {
					(0,-1.905199)(1,-1.905199)(2,-1.905199)(3,-1.905199)(4,-1.905199)(5, -1.905199) (6, -1.905199)(7, -1.905199)(8, -1.905199)(9, -1.905199)(10, -1.905199)  
				};
			\legend{CSL-tri, CSL-med, Centralize}
			\end{axis}
	\end{tikzpicture}}
	\subfigure{
		\begin{tikzpicture}
			\tikzstyle{every node}=[font=\tiny]
			\pgfplotsset{every axis legend/.append style={
					at={(0.5,1.03)},
					anchor=south},
				every axis y label/.append style={at={(-0.08,0.5)}}}
			\begin{axis}[title={small $n$, $t_{2.1}$ noises},
				xtick =data,xlabel = $T$,ylabel = {$\log(\text{error})$},legend columns=1,legend pos=north east,legend style={font=\footnotesize},font=\small,height=5.5cm, width=5.8cm]
				\addplot[
				color=blue,
				mark=triangle]
				coordinates
				{(0,-0.2561254)(1, -1.5151636)(2, -1.5523218)(3, -1.5546903 )(4,-1.5406540)(5, -1.5577291)(6,-1.5503611)(7, -1.5166572 )(8,-1.5357982)
				};
				\addplot[
				color=red,
				mark=square
				]
				coordinates {
					(0,-0.2490578)(1, -1.4318597)(2, -1.4941661)(3, -1.5422645 )(4,-1.5227891)(5, -1.5196542)(6,-1.5152028 )(7,-1.4909886)(8,-1.4926598)
				};
				\addplot[
				color=green,
				mark=o]
				coordinates  {(0,-1.905199)(1,-1.905199)(2,-1.905199)(3,-1.905199)(4,-1.905199)(5, -1.905199) (6, -1.905199)(7, -1.905199)(8, -1.905199)
				};
				\legend{CSL-tri, CSL-med, Centralize}
			\end{axis}
	\end{tikzpicture}}
\subfigure{
	\begin{tikzpicture}
		\tikzstyle{every node}=[font=\tiny]
		\pgfplotsset{every axis legend/.append style={
				at={(0.5,1.03)},
				anchor=south},
			every axis y label/.append style={at={(-0.08,0.5)}}}
		\begin{axis}[title={small $n$, $t_{2.1}$ noises},
			xtick =data,xlabel = $T$,ylabel = {$\log(\text{error})$},legend columns=1,legend pos=north east,legend style={font=\footnotesize},font=\small,height=5.5cm, width=5.8cm]
			\addplot[
			color=blue,
			mark=triangle]
			coordinates
			{(0,0.2360059)(1, -1.2402426)(2, -1.5192502)(3, -1.5291811)(4,  -1.5353815)(5,  -1.5372088)(6, -1.4639477)(7, -1.5126159)(8,  -1.5387646)(9, -1.4779335)
			};
			\addplot[
			color=red,
			mark=square
			]
			coordinates {
				(0,0.1622814 )(1, -1.2581293)(2, -1.4604702)(3, -1.4764791 )(4, -1.4756085)(5, -1.4818096)(6, -1.4279381)(7, -1.4595512)(8,  -1.4373875)(9,  -1.4488772)
			};
			\addplot[
			color=green,
			mark=o]
			coordinates  {(0,-1.905199)(1,-1.905199)(2,-1.905199)(3,-1.905199)(4,-1.905199)(5, -1.905199) (6, -1.905199)(7, -1.905199)(8, -1.905199)(9, -1.905199)
			};
			\legend{CSL-tri, CSL-med, Centralize}
		\end{axis}
\end{tikzpicture}}
\subfigure{
	\begin{tikzpicture}
		\tikzstyle{every node}=[font=\tiny]
		\pgfplotsset{every axis legend/.append style={
				at={(0.5,1.03)},
				anchor=south},
			every axis y label/.append style={at={(-0.08,0.5)}}}
		\begin{axis}[title={small $n$, $t_{2.1}$ noises},
			xtick =data,xlabel = $T$,ylabel = {$\log(\text{error})$},legend columns=1,legend pos=north east,legend style={font=\footnotesize},font=\small,height=5.5cm, width=5.8cm]
			\addplot[
			color=blue,
			mark=triangle]
			coordinates
			{(0,0.4203431)(1, -0.8361922)(2, -1.3213879)(3, -1.4428194)(4, -1.4659085)(5, -1.5074495)(6, -1.5087739)(7, -1.4668836)(8, -1.4508289 )(9,-1.4709829)(10,-1.4806844)
			};
			\addplot[
			color=red,
			mark=square
			]
			coordinates {
				(0,0.3906806)(1, -0.9903838)(2, -1.3172035)(3, -1.4191217 )(4,-1.3140600)(5, -1.3698999)(6,
				-1.3939758)(7, -1.4201573)(8, -1.3899591)
				(9, -1.4153187 )(10,-1.4522026)
			};
			\addplot[
			color=green,
			mark=o]
			coordinates  {(0,-1.905199)(1,-1.905199)(2,-1.905199)(3,-1.905199)(4,-1.905199)(5, -1.905199) (6, -1.905199)(7, -1.905199)(8, -1.905199)(9, -1.905199)(10,-1.905199)
			};
			\legend{CSL-tri, CSL-med, Centralize}
		\end{axis}
\end{tikzpicture}}
\caption{Estimation error $\log\|\widehat{\boldsymbol{\Theta}}_{t}-\boldsymbol{\Theta}^{\star}\|_{F}$ versus iteration count $T$ under varying local sample sizes, noise distributions and Byzantine fractions. Three curves in each panel correspond to CSL-tri (Option~I), CSL-med (Option~II) and the centralized $\ell_{1}$-penalized oracle.}\label{fig2}
\end{figure}
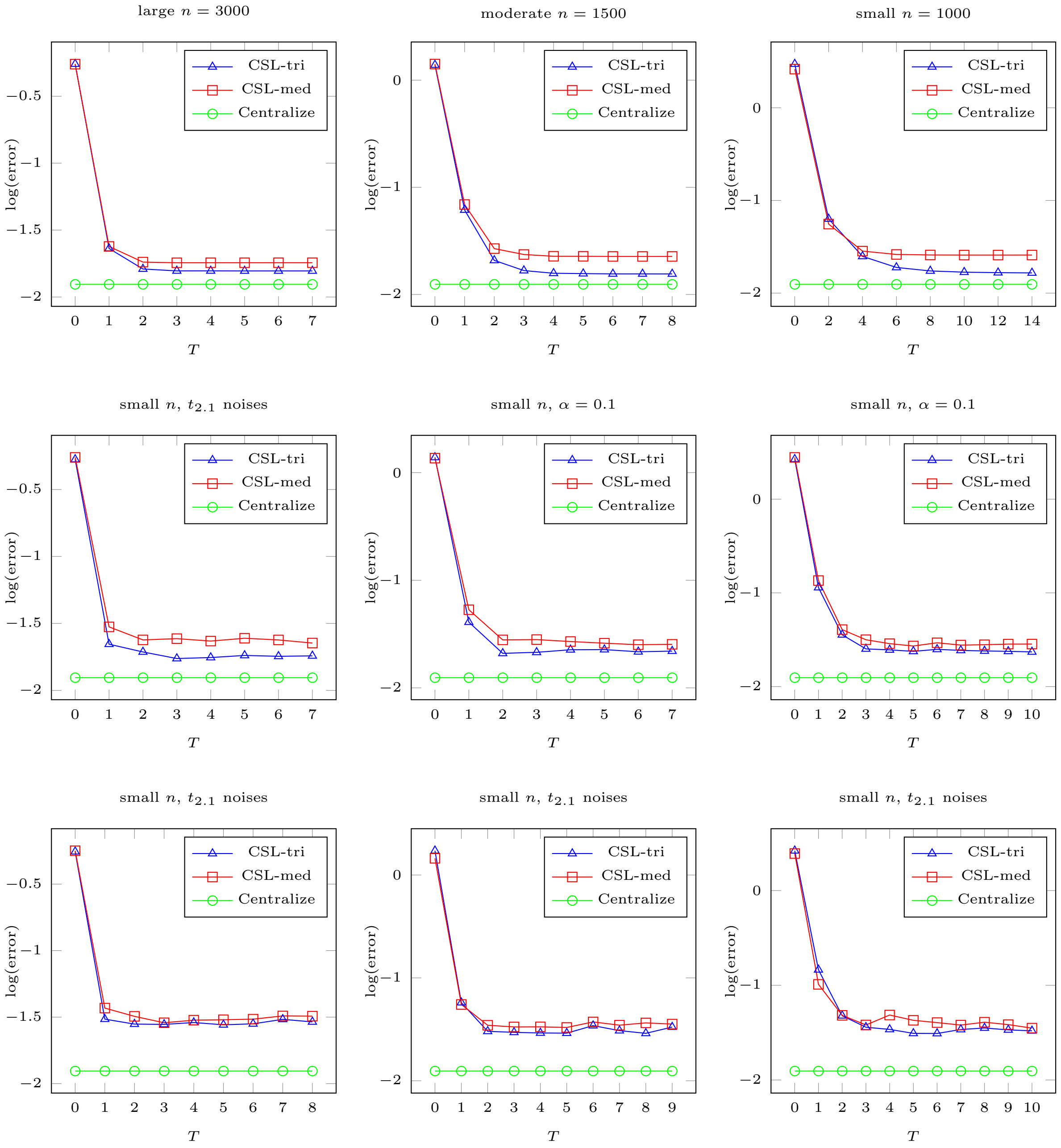

We next turn to the behaviour of Algorithm~\ref{Al2} under varying local sample sizes and Byzantine attack rates. The two robust aggregation options are labelled \textbf{CSL-tri} (Option~I, trimmed mean) and \textbf{CSL-med} (Option~II, coordinate-wise median); the centralized $\ell_{1}$-penalized MLE based on the full sample serves as an oracle baseline (\textbf{Centralize}). Figure~\ref{fig2} demonstrates the log estimation error $\log\|\widehat{\boldsymbol{\Theta}}_{t}-\boldsymbol{\Theta}^{\star}\|_{F}$ against the iteration count $T$ across nine configurations of local sample size $n$, noise distribution and Byzantine fraction $\alpha$. We set $K=4$, $s=5$, $d=50$, total sample size $N=30000$, and regularization parameter $\lambda=4\sqrt{\log(d)/N}+2\sqrt{\log(d)/n},4\sqrt{\log(d)/N}+3\sqrt{\log(d)/n}$.

\paragraph{Byzantine attack setting.}
In these experiments, Byzantine workers corrupt their local gradient propagates before aggregation. The simulations vary the local sample size, noise distribution, and Byzantine fraction, demonstrating both optimization convergence and the effect of contamination. CSL-tri and CSL-med are evaluated under the same attack realization.

Corrupted gradients from Byzantine machines are generated by three representative attack modes:
\begin{itemize}
\item[-] Gaussian attack. The true gradient is replaced by a random vector whose entries are independent draws from $N(0,10)$.
\item[-] Omniscient attack. The true gradient is replaced by its sign-flipped negative, modelling an adversary that knows the gradient and pushes the iterate in the opposite direction.
For Byzantine machine $j$, assume the true gradient is $\nabla \hat{L}_j$, machine $j$ reports $-5 \nabla \hat{L}_j$.
\item[-] Bit-flip attack. The true gradient is replaced by its coordinate-wise sign-flipped copy.
For Byzantine machine $j$, assume the true gradient is $\nabla \hat{L}_j$. Machine $j$ takes its first ten coordinates multiplied with -10 and reports to the central machine.
\end{itemize}
The Byzantine fraction is varied over $\alpha \in\{0.1,0.2\}$.

Figure~\ref{fig2} first isolates the effect of the local sample size and the attack/noise configuration on the outer iterations. The error decreases sharply during the initial updates and then enters a plateau, with larger local samples generally producing a lower starting error and a tighter final band. The two robust curves remain close to one another across the panels, while the centralized reference is constant because it is fitted once. The panels with heavy-tailed noise or a nonzero Byzantine fraction show a larger residual gap and more visible fluctuations, but the robust trajectories remain bounded and retain the same early-improvement/late-stabilization pattern.

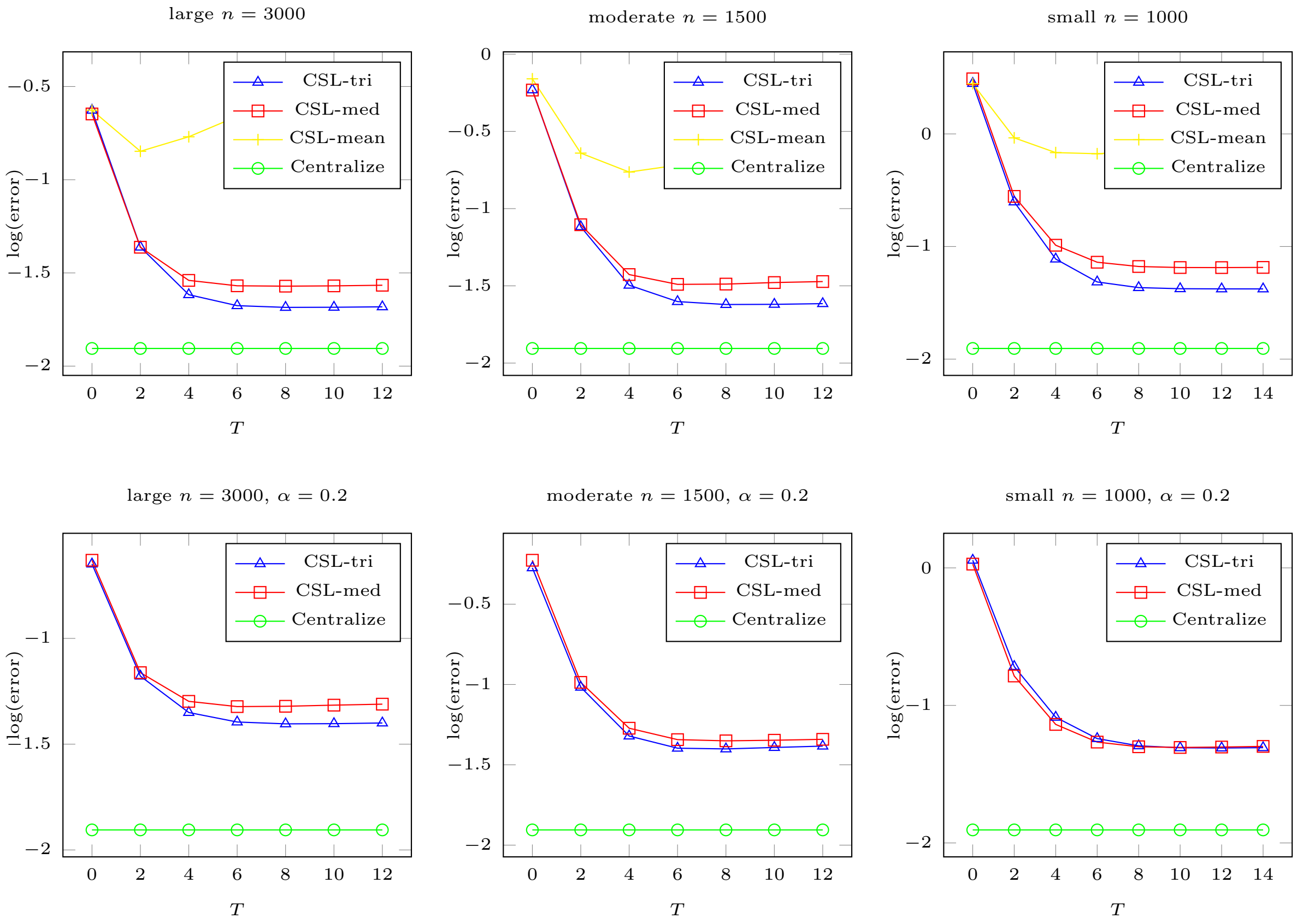
\begin{figure}
	\centering
	\subfigure{
		\begin{tikzpicture}
			\tikzstyle{every node}=[font=\tiny]
			\pgfplotsset{every axis legend/.append style={
					at={(0.5,1.03)},
					anchor=south},
				every axis y label/.append style={at={(-0.08,0.5)}}}
			\begin{axis}[title={large $n=3000$},  
				xtick =data, xlabel = $T$,ylabel = {$\log(\text{error})$}, legend columns=1,legend pos=north east,legend style={font=\footnotesize},font=\small,height=5.5cm, width=5.8cm]
				\addplot[
				color=blue,
				mark=triangle
				] 
				coordinates {
					(0,-0.6273444)(2, -1.3615201)(4,-1.6167562)(6,-1.6756178)(8, -1.6848851)(10,-1.6841802)(12,-1.6818118) };
				\addplot[
				color=red,
				mark=square
				]
				coordinates {
					(0,-0.6471292)(2, -1.3621783)(4, -1.5400087)(6, -1.5690242)(8,-1.5710233)(10, -1.5694049) (12,  -1.5660569) };
				\addplot[
			color=yellow,
			mark=+
			]
			coordinates {
				(0,-0.6268152)(2,-0.8474323)(4, -0.7695140)(6, -0.6603902)(8,-0.5730577)(10, -0.5076838) (12,-0.4586218)
			};
                \addplot[
				color=green,
				mark=o
				]
				coordinates {
				(0,-1.905199)(2, -1.905199)(4, -1.905199)(6, -1.905199)(8,-1.905199)(10, -1.905199) (12, -1.905199)
				};
				\legend{CSL-tri, CSL-med, CSL-mean, Centralize}
			\end{axis}
	\end{tikzpicture}}
	\subfigure{
		\begin{tikzpicture}
			\tikzstyle{every node}=[font=\tiny]
			\pgfplotsset{every axis legend/.append style={
					at={(0.5,1.03)},
					anchor=south},
				every axis y label/.append style={at={(-0.08,0.5)}}}
			\begin{axis}[title={moderate $n=1500$},
				xtick =data,xlabel = $T$,ylabel = {$\log(\text{error})$},legend columns=1,legend pos=north east,legend style={font=\footnotesize},font=\small,height=5.5cm, width=5.8cm]
				\addplot[
				color=blue,
				mark=triangle]
				coordinates
				{(0,-0.2314149)(2, -1.1172129)(4,-1.4957910)(6, -1.6018347)(8,-1.6207718)(10, -1.6195988)(12, -1.6150712)}; 
				\addplot[
				color=red,
				mark=square
				]
				coordinates {
					(0,-0.2325288)(2, -1.1042839)(4, -1.4268728)(6,  -1.4903684)(8,-1.4885167)(10, -1.4784970)(12,-1.4722234)};
             \addplot[
			color=yellow,
			mark=+
			]
			coordinates {
				(0,-0.1596806)(2,-0.6407021)(4, -0.7631552)(6, -0.7187355)(8,-0.6465574)(10, -0.5865886) (12,-0.5411773)
			};
				\addplot[
				color=green,
				mark=o
				]
				coordinates {
					(0,-1.905199)(2, -1.905199)(4, -1.905199)(6, -1.905199)(8,-1.905199)(10, -1.905199) (12, -1.905199)
				};
				\legend{CSL-tri, CSL-med,CSL-mean, Centralize}
			\end{axis}
	\end{tikzpicture}}
	\subfigure{
		\begin{tikzpicture}
			\tikzstyle{every node}=[font=\tiny]
			\pgfplotsset{every axis legend/.append style={
					at={(0.5,1.03)},
					anchor=south},
				every axis y label/.append style={at={(-0.08,0.5)}}}
			\begin{axis}[title={small $n=1000$},
				xtick =data,xlabel = $T$,ylabel = {$\log(\text{error})$},legend columns=1,legend pos=north east,legend style={font=\footnotesize},font=\small,height=5.5cm, width=5.8cm]
				\addplot[
				color=blue,
				mark=triangle
				]
				coordinates {
					(0,0.4481351)(2,-0.6041786)(4,-1.1112232)(6,-1.3146614)(8, -1.3649024)(10, -1.3747474)(12,-1.3762893)(14,-1.3764633)
				};
				\addplot[
				color=red,
				mark=square
				]
				coordinates {(0,0.4892295)(2,-0.5550022)(4,-0.9886074)(6,-1.1395316)(8, -1.1775962)(10, -1.1866387)(12,-1.1871381)(14,-1.1855636)};
			\addplot[
			color=yellow,
			mark=+
			]
			coordinates {
				(0,0.44594167)(2,-0.03503341)(4, -0.16604911)(6, -0.17613223)(8,-0.15468129)(10, -0.12911082) (12,-0.10636531)(14,-0.08751106)
			};
            \addplot[
			color=green,
			mark=o
			]
			coordinates {
				(0,-1.905199)(2, -1.905199)(4, -1.905199)(6, -1.905199)(8,-1.905199)(10, -1.905199) (12, -1.905199)(14,-1.905199)
			};
				\legend{CSL-tri, CSL-med, CSL-mean,Centralize}
			\end{axis}
	\end{tikzpicture}}
	\subfigure{
		\begin{tikzpicture}
			\tikzstyle{every node}=[font=\tiny]
			\pgfplotsset{every axis legend/.append style={
					at={(0.5,1.03)},
					anchor=south},
				every axis y label/.append style={at={(-0.08,0.5)}}}
			\begin{axis}[title={large $n=3000$, $\alpha=0.2$},
				xtick =data,xlabel = $T$,ylabel = {$\log(\text{error})$},legend columns=1,legend pos=north east,legend style={font=\footnotesize},font=\small,height=5.5cm, width=5.8cm]
                \addplot[
				color=blue,
				mark=triangle
				]
				coordinates {
					(0,-0.6484125)(2,-1.1782506)(4,-1.3502431)(6,-1.3951776)(8,-1.4038924)(10, -1.4030149)(12,-1.3998298)
				};
				\addplot[
				color=red,
				mark=square
				]
				coordinates {
					(0, -0.6306325)(2, -1.1619361)(4, -1.2975230)(6, -1.3223996)(8,-1.3208057)(10, -1.3154932)(12,  -1.3108346)
				};
				\addplot[
				color=green,
				mark=o
				]
				coordinates {
					(0,-1.905199)(2, -1.905199)(4, -1.905199)(6, -1.905199)(8,-1.905199)(10, -1.905199) (12, -1.905199)
				};
				\legend{CSL-tri, CSL-med, Centralize}
			\end{axis}
	\end{tikzpicture}}
\subfigure{
	\begin{tikzpicture}
		\tikzstyle{every node}=[font=\tiny]
		\pgfplotsset{every axis legend/.append style={
				at={(0.5,1.03)},
				anchor=south},
			every axis y label/.append style={at={(-0.08,0.5)}}}
		\begin{axis}[title={moderate $n=1500$, $\alpha=0.2$},
			xtick =data,xlabel = $T$,ylabel = {$\log(\text{error})$},legend columns=1,legend pos=north east,legend style={font=\footnotesize},font=\small,height=5.5cm, width=5.8cm]
			\addplot[
			color=blue,
			mark=triangle] 
			coordinates {(0, -0.2735128)(2, -1.0168565)(4, -1.3208024)(6, -1.3965508)(8, -1.4010402)(10, -1.3920353)(12, -1.3835264)
			};
			\addplot[
			color=red,
			mark=square
			]
			coordinates {
				(0,-0.2265322)(2, -0.9869219)(4, -1.2721642)(6, -1.3432349)(8,-1.3507279)(10, -1.3470459)(12, -1.3413283)
			};
			\addplot[
			color=green,
			mark=o
			]
			coordinates {
					(0,-1.905199)(2, -1.905199)(4, -1.905199)(6, -1.905199)(8,-1.905199)(10, -1.905199) (12, -1.905199)
			};
		\legend{CSL-tri, CSL-med, Centralize}
		\end{axis}
\end{tikzpicture}}
	\subfigure{
		\begin{tikzpicture}
			\tikzstyle{every node}=[font=\tiny]
			\pgfplotsset{every axis legend/.append style={
					at={(0.5,1.03)},
					anchor=south},
				every axis y label/.append style={at={(-0.08,0.5)}}}
			\begin{axis}[title={small $n=1000$, $\alpha=0.2$},
				xtick =data,xlabel = $T$,ylabel = {$\log(\text{error})$},legend columns=1,legend pos=north east,legend style={font=\footnotesize},font=\small,height=5.5cm, width=5.8cm]
				\addplot[
				color=blue,
				mark=triangle
				]
				coordinates {
					(0,0.05545051)(2,-0.71727228)(4, -1.08618733)(6, -1.24118480 )(8,-1.29373502)(10, -1.30799517)(12,-1.30970192)(14, -1.30729125)
				};
				\addplot[
				color=red,
				mark=square
				]
				coordinates {
					(0, 0.02753279)(2, -0.78601432)(4, -1.13814807)(6, -1.26669783)(8,-1.30061056)(10, -1.30571669)(12, -1.30270746)(14,-1.29802859)
				};
				\addplot[
				color=green,
				mark=o
				]
				coordinates {
					(0,-1.905199)(2, -1.905199)(4, -1.905199)(6, -1.905199)(8,-1.905199)(10, -1.905199) (12, -1.905199)(14,-1.905199)
				};
			\legend{CSL-tri, CSL-med, Centralize}
			\end{axis}
	\end{tikzpicture}}
\caption{Estimation error $\log\|\widehat{\boldsymbol{\Theta}}_{t}-\boldsymbol{\Theta}^{\star}\|_{F}$ versus iteration count $T$ under varying local sample sizes, noise distributions and Byzantine fractions. Three curves in each panel correspond to CSL-tri (Option~I), CSL-med (Option~II) and the centralized $\ell_{1}$-penalized oracle.}\label{fig3}
\end{figure}

Figure~\ref{fig3} shows the same comparison under the alternative configurations in the middle group of panels. Relative to the cleaner, larger-sample cases, the robust methods typically start farther from the oracle and converge more slowly, especially when the local sample size is small. Nevertheless, both aggregators continue to reduce the error over the outer iterations rather than being driven away by the contaminated reports; the remaining difference between trimming and the median is secondary to this common stability behavior.

\begin{figure}
	\centering
	\subfigure{
		\begin{tikzpicture}
			\tikzstyle{every node}=[font=\tiny]
			\pgfplotsset{every axis legend/.append style={
					at={(0.5,1.03)},
					anchor=south},
				every axis y label/.append style={at={(-0.08,0.5)}}}
			\begin{axis}[title={large $n=3000$},  
				xtick =data, xlabel = $T$,ylabel = {$\log(\text{error})$}, legend columns=1,legend pos=north east,legend style={font=\footnotesize},font=\small,height=5.5cm, width=5.8cm]
				\addplot[
				color=blue,
				mark=triangle
				]
				coordinates {
					(0,-0.6468722)(2, -1.4753949)(4,-1.7831925)(6, -1.8644112)(8,-1.8824035)(10, -1.8842915) (12,  -1.8828827)  };
				\addplot[
				color=red,
				mark=square
				] 
				coordinates {
					(0,-0.6361188)(2, -1.4436155)(4,-1.6483782)(6, -1.6862996)(8,-1.6910845)(10, -1.6909977) (12,  -1.6905139)  };
				\addplot[
				color=green,
				mark=o
				]
				coordinates {
				(0,-1.905199)(2, -1.905199)(4, -1.905199)(6, -1.905199)(8,-1.905199)(10, -1.905199) (12, -1.905199)
				};
				
				\legend{CSL-tri, CSL-med, Centralize}
			\end{axis}
	\end{tikzpicture}}
	\subfigure{
		\begin{tikzpicture}
			\tikzstyle{every node}=[font=\tiny]
			\pgfplotsset{every axis legend/.append style={
					at={(0.5,1.03)},
					anchor=south},
				every axis y label/.append style={at={(-0.08,0.5)}}}
			\begin{axis}[title={moderate $n=1500$},
				xtick =data,xlabel = $T$,ylabel = {$\log(\text{error})$},legend columns=1,legend pos=north east,legend style={font=\footnotesize},font=\small,height=5.5cm, width=5.8cm]
				\addplot[
				color=blue,
				mark=triangle]
				coordinates
				{(0,-0.2445000)(2, -1.2083493)(4,-1.6307850)(6, -1.7467096)(8,-1.7662104)(10, -1.7689031)(12, -1.7694013)}; 
				\addplot[
				color=red,
				mark=square
				]
				coordinates {
					(0,-0.2280503)(2, -1.2094258)(4,-1.5963044)(6, -1.6891696)(8,-1.7002263)(10, -1.6961455)(12, -1.6920257) };
				\addplot[
				color=green,
				mark=o
				]
				coordinates {
					(0,-1.905199)(2, -1.905199)(4, -1.905199)(6, -1.905199)(8,-1.905199)(10, -1.905199) (12, -1.905199)
				};
				\legend{CSL-tri, CSL-med, Centralize}
			\end{axis}
	\end{tikzpicture}}
	\subfigure{
		\begin{tikzpicture}
			\tikzstyle{every node}=[font=\tiny]
			\pgfplotsset{every axis legend/.append style={
					at={(0.5,1.03)},
					anchor=south},
				every axis y label/.append style={at={(-0.08,0.5)}}}
			\begin{axis}[title={small $n=1000$},
				xtick =data,xlabel = $T$,ylabel = {$\log(\text{error})$},legend columns=1,legend pos=north east,legend style={font=\footnotesize},font=\small,height=5.5cm, width=5.8cm]
				\addplot[
				color=blue,
				mark=triangle
				]
				coordinates {
					(0,0.0452620)(2,-1.0202014)(4,-1.5585811)(6,-1.6985082)(8, -1.7058575)(10, -1.6968175)(12,-1.6901495)(14,-1.6860926)
				};
				\addplot[
				color=red,
				mark=square
				]
				coordinates {(0,0.02846986)(2,-0.98399606)(4,-1.42816031)(6,-1.56147793)(8, -1.58936334)(10, -1.59553717)(12,-1.59715027)(14,-1.59716863)};
            \addplot[
			color=green,
			mark=o
			]
			coordinates {
			(0,-1.905199)(2, -1.905199)(4, -1.905199)(6, -1.905199)(8,-1.905199)(10, -1.905199) (12, -1.905199)(14,-1.905199)
			};
				\legend{CSL-tri, CSL-med, Centralize}
			\end{axis}
	\end{tikzpicture}}
	\subfigure{
		\begin{tikzpicture}
			\tikzstyle{every node}=[font=\tiny]
			\pgfplotsset{every axis legend/.append style={
					at={(0.5,1.03)},
					anchor=south},
				every axis y label/.append style={at={(-0.08,0.5)}}}
			\begin{axis}[title={small $n$, $t_{2.1}$ noises},
				xtick =data,xlabel = $T$,ylabel = {$\log(\text{error})$},legend columns=1,legend pos=north east,legend style={font=\footnotesize},font=\small,height=5.5cm, width=5.8cm]
				\addplot[
				color=blue,
				mark=triangle]
				coordinates {(0,-0.5911667)(2,-1.3670906)(4,-1.6192278)(6,-1.6812300)(8, -1.6934598)(10, -1.6942528)(12,-1.6926634)
				};
				\addplot[
				color=red,
				mark=square
				]
				coordinates {
					(0, -0.6294708)(2,-1.3726278)(4,-1.5711754)(6,-1.6101630)(8,-1.6149878)(10, -1.6129532)(12,-1.6103375)
				};
				\addplot[
				color=green,
				mark=o
				]
				coordinates {
					(0,-1.905199)(2, -1.905199)(4, -1.905199)(6, -1.905199)(8,-1.905199)(10, -1.905199) (12, -1.905199)
				};
				\legend{CSL-tri, CSL-med, Centralize}
			\end{axis}
	\end{tikzpicture}}
\subfigure{
	\begin{tikzpicture}
		\tikzstyle{every node}=[font=\tiny]
		\pgfplotsset{every axis legend/.append style={
				at={(0.5,1.03)},
				anchor=south},
			every axis y label/.append style={at={(-0.08,0.5)}}}
		\begin{axis}[title={small $n$, $\alpha=0.1$},
			xtick =data,xlabel = $T$,ylabel = {$\log(\text{error})$},legend columns=1,legend pos=north east,legend style={font=\footnotesize},font=\small,height=5.5cm, width=5.8cm]
			\addplot[
			color=blue,
			mark=triangle]
			coordinates {(0, -0.2400974)(2, -1.1932976)(4, -1.5757020)(6, -1.6593969)(8, -1.6684869)(10, -1.6650242)(12, -1.6611587)
			};
			\addplot[
			color=red,
			mark=square
			] 
			coordinates {
				(0,-0.2197867)(2, -1.1294530)(4, -1.4857705)(6, -1.5817456)(8,-1.5989733)(10, -1.5996417)(12, -1.5975411)
			};
			\addplot[
			color=green,
			mark=o
			]
			coordinates {
			(0,-1.905199)(2, -1.905199)(4, -1.905199)(6, -1.905199)(8,-1.905199)(10, -1.905199) (12, -1.905199)
			};
		\legend{CSL-tri, CSL-med, Centralize}
		\end{axis}
\end{tikzpicture}}
	\subfigure{
		\begin{tikzpicture}
			\tikzstyle{every node}=[font=\tiny]
			\pgfplotsset{every axis legend/.append style={
					at={(0.5,1.03)},
					anchor=south},
				every axis y label/.append style={at={(-0.08,0.5)}}}
			\begin{axis}[title={small $n$, $\alpha=0.1$},
				xtick =data,xlabel = $T$,ylabel = {$\log(\text{error})$},legend columns=1,legend pos=north east,legend style={font=\footnotesize},font=\small,height=5.5cm, width=5.8cm]
				\addplot[
				color=blue,
				mark=triangle
				]
				coordinates {
					(0,0.03504988)(2,-0.97820571)(4, -1.47109942)(6, -1.62456644)(8,-1.64819388)(10, -1.64624191)(12,-1.64095250)
				};
				\addplot[
				color=red,
				mark=square
				]
				coordinates {
					(0, 0.07156963)(2, -0.95548489)(4, -1.43183370)(6, -1.56975210)(8,-1.58986724)(10, -1.58758599
					)(12, -1.58115983)
				};
				\addplot[
				color=green,
				mark=o
				]
				coordinates {
					(0,-1.905199)(2, -1.905199)(4, -1.905199)(6, -1.905199)(8,-1.905199)(10, -1.905199) (12, -1.905199)
				};
			\legend{CSL-tri, CSL-med, Centralize}
			\end{axis}
	\end{tikzpicture}}
\caption{Estimation error $\log\|\widehat{\boldsymbol{\Theta}}_{t}-\boldsymbol{\Theta}^{\star}\|_{F}$ versus iteration count $T$ under varying local sample sizes, noise distributions and Byzantine fractions. Three curves in each panel correspond to CSL-tri (Option~I), CSL-med (Option~II) and the centralized $\ell_{1}$-penalized oracle.}\label{fig4}
\end{figure}
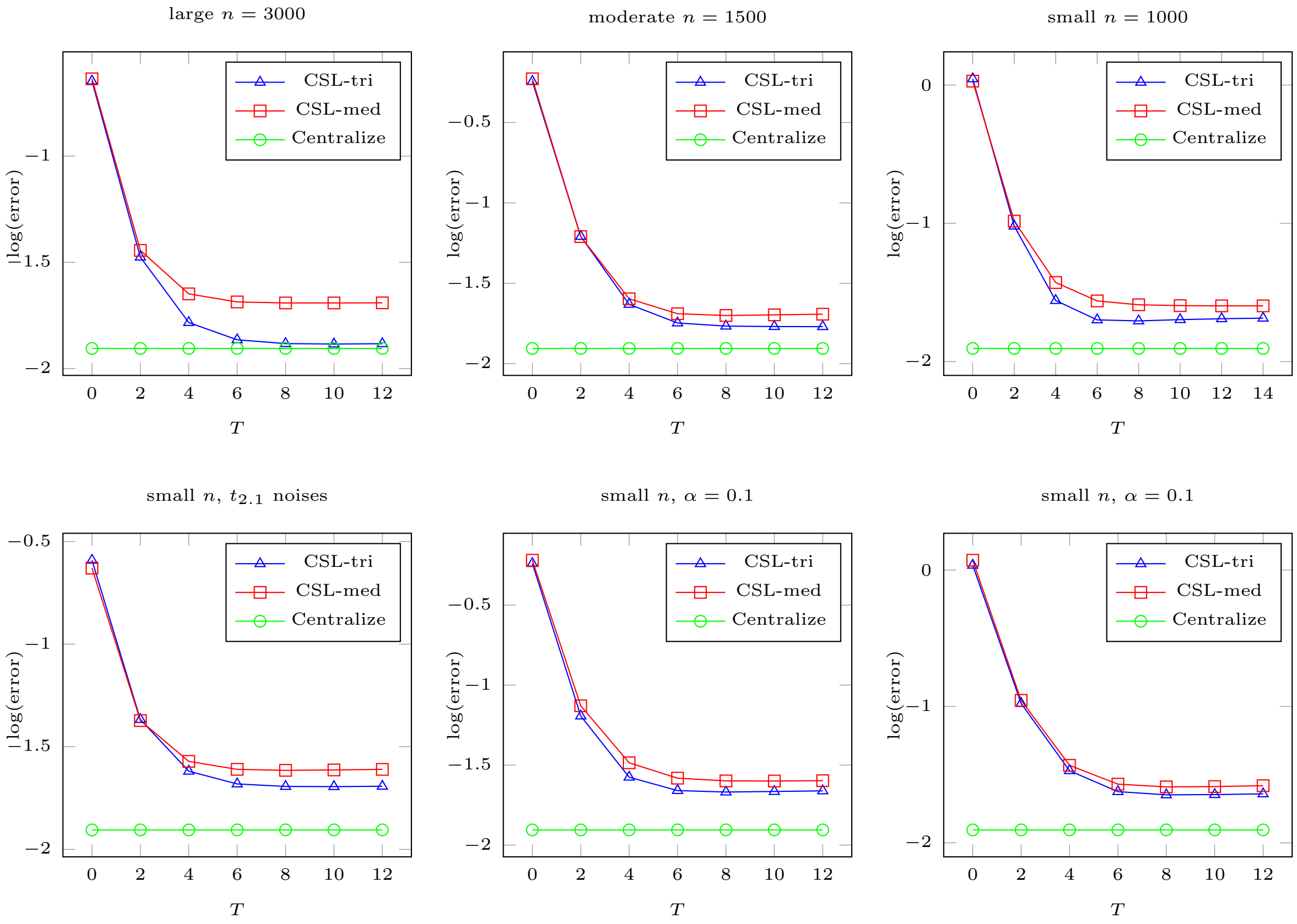
The final group in Figure~\ref{fig4} emphasizes the more difficult small-sample, contaminated settings. Here the curves may exhibit a mild late-iteration fluctuation instead of a perfectly monotone decrease, reflecting the higher variance of local gradients. The important feature is that the fluctuation remains around a stable error level: neither robust rule shows the unbounded drift that would result if a small number of adversarial reports controlled the aggregate. Taken together, Figures~\ref{fig2}--\ref{fig4} show a consistent pattern of fast initial improvement, a contamination- and sample-size-dependent residual error, and bounded long-run behavior.


\subsection{Experiments on Real Data}\label{sec:exp_real}
We evaluate the proposed method on the Connect-4 data set distributed through the LIBSVM classification repository. The data originate from the UCI Connect-4 data set, encodeing the $42$ board positions using three binary indicators for the symbols $\{o,b,x\}$, resulting in $67{,}557$ observations, $3$ classes, and $126$ sparse features.\footnote{See the LIBSVM multiclass data repository: \url{https://www.csie.ntu.edu.tw/~cjlin/libsvmtools/datasets/multiclass.html}.} This setting combines a moderate sample size with a high-dimensional, structured feature representation.

Following the methods and notation introduced above, we compare the robust distributed procedures with the non-robust and centralized references. The goal here is to test whether the trajectory behavior observed in the controlled experiments persists on a real, sparse, multiclass data set. Accuracy and multinomial logloss are reported over the outer iterations and at the final iterate; results are averaged over $20$ repeats with five outer folds per repeat. Runtime measures optimization only and excludes validation and test evaluation.

\begin{figure}[!ht]
  \centering
  \subfigure{
    \begin{tikzpicture}
      \tikzstyle{every node}=[font=\tiny]
      \pgfplotsset{every axis legend/.append style={at={(0.5,1.20)},anchor=south}, every axis y label/.append style={at={(-0.08,0.5)}}}
      \begin{axis}[title={$\alpha=0.2$}, xtick={0,4,8,12,16,20}, xlabel={$T$}, ylabel={Accuracy}, ymin=0.64, ymax=0.76, legend columns=1, legend style={font=\footnotesize}, legend pos=south east, font=\small, height=5.5cm, width=7.9cm, mark repeat=4]
        \addplot[forget plot,color=blue,mark=triangle,solid] coordinates {
          (0,0.66767546)(1,0.67757452)(2,0.68959032)(3,0.70186144)(4,0.71252280)(5,0.72111924)(6,0.72776028)(7,0.73269243)(8,0.73662172)(9,0.73968580)(10,0.74227178)(11,0.74411171)(12,0.74562821)(13,0.74687973)(14,0.74796474)(15,0.74881956)(16,0.74950120)(17,0.75011475)(18,0.75059286)(19,0.75086669)(20,0.75115533)
        };
        \addplot[forget plot,color=red,mark=square,solid] coordinates {
          (0,0.66767546)(1,0.67722593)(2,0.68861190)(3,0.70047004)(4,0.71109883)(5,0.71977004)(6,0.72673380)(7,0.73214035)(8,0.73629018)(9,0.73951560)(10,0.74201794)(11,0.74393780)(12,0.74549722)(13,0.74674727)(14,0.74761764)(15,0.74845693)(16,0.74909860)(17,0.74970845)(18,0.75014067)(19,0.75039749)(20,0.75067872)
        };
        \addplot[forget plot,color=green,mark=o,solid] coordinates {
          (0,0.74488357)(1,0.74488357)(2,0.74488357)(3,0.74488357)(4,0.74488357)(5,0.74488357)(6,0.74488357)(7,0.74488357)(8,0.74488357)(9,0.74488357)(10,0.74488357)(11,0.74488357)(12,0.74488357)(13,0.74488357)(14,0.74488357)(15,0.74488357)(16,0.74488357)(17,0.74488357)(18,0.74488357)(19,0.74488357)(20,0.74488357)
        };
        \addplot[forget plot,color=blue,mark=x,dashed,mark options={solid}] coordinates {
          (0,0.65886584)(1,0.66367509)(2,0.67421505)(3,0.68836903)(4,0.70174668)(5,0.71258052)(6,0.72105781)(7,0.72739172)(8,0.73212327)(9,0.73557591)(10,0.73810266)(11,0.74009063)(12,0.74150351)(13,0.74268030)(14,0.74373644)(15,0.74457942)(16,0.74526698)(17,0.74579764)(18,0.74625207)(19,0.74667614)(20,0.74704546)
        };
        \addplot[forget plot,color=red,mark=+,dashed,mark options={solid}] coordinates {
          (0,0.65886584)(1,0.66326729)(2,0.67242472)(3,0.68518874)(4,0.69782696)(5,0.70855422)(6,0.71742746)(7,0.72426686)(8,0.72960163)(9,0.73356052)(10,0.73655061)(11,0.73897819)(12,0.74087290)(13,0.74241827)(14,0.74365057)(15,0.74451948)(16,0.74529365)(17,0.74593459)(18,0.74627579)(19,0.74659256)(20,0.74682497)
        };
        \addlegendimage{blue,mark=triangle,solid}
        \addlegendentry{CSL-tri, $m=10$}
        \addlegendimage{blue,mark=x,dashed,mark options={solid}}
        \addlegendentry{CSL-tri, $m=20$}
        \addlegendimage{red,mark=square,solid}
        \addlegendentry{CSL-med, $m=10$}
        \addlegendimage{red,mark=+,dashed,mark options={solid}}
        \addlegendentry{CSL-med, $m=20$}
        \addlegendimage{green,mark=o,solid}
        \addlegendentry{Centralize}
      \end{axis}
    \end{tikzpicture}}
  \hspace{0.2cm}%
  \subfigure{
    \begin{tikzpicture}
      \tikzstyle{every node}=[font=\tiny]
      \pgfplotsset{every axis legend/.append style={at={(0.5,1.20)},anchor=south}, every axis y label/.append style={at={(-0.08,0.5)}}}
      \begin{axis}[title={$\alpha=0.3$}, xtick={0,4,8,12,16,20}, xlabel={$T$}, ylabel={Accuracy}, ymin=0.64, ymax=0.76, legend columns=1, legend style={font=\footnotesize}, legend pos=south east, font=\small, height=5.5cm, width=7.9cm, mark repeat=4]
        \addplot[forget plot,color=blue,mark=triangle,solid] coordinates {
          (0,0.66767546)(1,0.67678037)(2,0.68788803)(3,0.69867076)(4,0.70850243)(5,0.71633507)(6,0.72268898)(7,0.72769810)(8,0.73185458)(9,0.73512070)(10,0.73783618)(11,0.74007058)(12,0.74196157)(13,0.74338703)(14,0.74466892)(15,0.74570138)(16,0.74650145)(17,0.74721419)(18,0.74772487)(19,0.74818892)(20,0.74861153)
        };
        \addplot[forget plot,color=red,mark=square,solid] coordinates {
          (0,0.66767546)(1,0.67686177)(2,0.68783029)(3,0.69884398)(4,0.70881187)(5,0.71729656)(6,0.72394279)(7,0.72908954)(8,0.73329043)(9,0.73681042)(10,0.73936011)(11,0.74152421)(12,0.74326569)(13,0.74455053)(14,0.74558300)(15,0.74639343)(16,0.74700551)(17,0.74749770)(18,0.74787738)(19,0.74830591)(20,0.74868781)
        };
        \addplot[forget plot,color=green,mark=o,solid] coordinates {
          (0,0.74488357)(1,0.74488357)(2,0.74488357)(3,0.74488357)(4,0.74488357)(5,0.74488357)(6,0.74488357)(7,0.74488357)(8,0.74488357)(9,0.74488357)(10,0.74488357)(11,0.74488357)(12,0.74488357)(13,0.74488357)(14,0.74488357)(15,0.74488357)(16,0.74488357)(17,0.74488357)(18,0.74488357)(19,0.74488357)(20,0.74488357)
        };
        \addplot[forget plot,color=blue,mark=x,dashed,mark options={solid}] coordinates {
          (0,0.65886584)(1,0.66318957)(2,0.67247504)(3,0.68447528)(4,0.69621945)(5,0.70623399)(6,0.71455215)(7,0.72112511)(8,0.72598693)(9,0.73021225)(10,0.73328743)(11,0.73584822)(12,0.73788723)(13,0.73951326)(14,0.74080772)(15,0.74184389)(16,0.74257587)(17,0.74317686)(18,0.74382594)(19,0.74425892)(20,0.74461271)
        };
        \addplot[forget plot,color=red,mark=+,dashed,mark options={solid}] coordinates {
          (0,0.65886584)(1,0.66310076)(2,0.67164315)(3,0.68312971)(4,0.69518471)(5,0.70564257)(6,0.71405771)(7,0.72048342)(8,0.72547994)(9,0.72960979)(10,0.73282707)(11,0.73529165)(12,0.73751273)(13,0.73930751)(14,0.74049317)(15,0.74150269)(16,0.74221763)(17,0.74293629)(18,0.74356539)(19,0.74411309)(20,0.74461638)
        };
        \addlegendimage{blue,mark=triangle,solid}
        \addlegendentry{CSL-tri, $m=10$}
        \addlegendimage{blue,mark=x,dashed,mark options={solid}}
        \addlegendentry{CSL-tri, $m=20$}
        \addlegendimage{red,mark=square,solid}
        \addlegendentry{CSL-med, $m=10$}
        \addlegendimage{red,mark=+,dashed,mark options={solid}}
        \addlegendentry{CSL-med, $m=20$}
        \addlegendimage{green,mark=o,solid}
        \addlegendentry{Centralize}
      \end{axis}
    \end{tikzpicture}}
  \caption{Test Accuracy trajectories on Connect-4 under Byzantine fractions $\alpha=0.2$ and $\alpha=0.3$. Solid curves denote $m=10$ and dashed curves denote $m=20$. CSL-tri is the coordinate-wise trimmed mean and CSL-med is the coordinate-wise median.}
  \label{fig:real-connect4-accuracy-step4}
\end{figure}
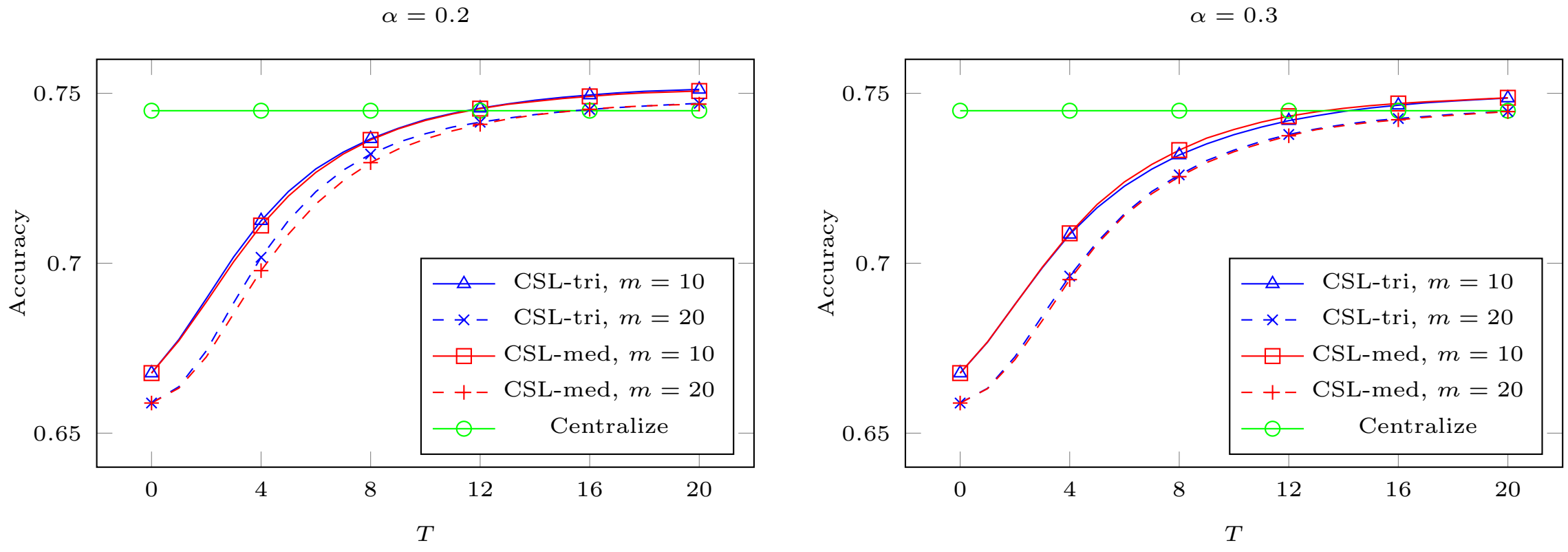

\begin{figure}[!ht]
  \centering
  \subfigure{
    \begin{tikzpicture}
      \tikzstyle{every node}=[font=\tiny]
      \pgfplotsset{every axis legend/.append style={at={(0.5,1.20)},anchor=south}, every axis y label/.append style={at={(-0.08,0.5)}}}
      \begin{axis}[title={$\alpha=0.2$}, xtick={0,4,8,12,16,20}, xlabel={$T$}, ylabel={Logloss}, ymin=0.6, ymax=0.84, legend columns=1, legend style={font=\footnotesize}, legend pos=north east, font=\small, height=5.5cm, width=7.9cm, mark repeat=4]
        \addplot[forget plot,color=blue,mark=triangle,solid] coordinates {
          (0,0.77163194)(1,0.74349703)(2,0.72299648)(3,0.70719966)(4,0.69471939)(5,0.68364910)(6,0.67510619)(7,0.66752940)(8,0.66096845)(9,0.65544049)(10,0.65053156)(11,0.64646269)(12,0.64314183)(13,0.63996096)(14,0.63789596)(15,0.63533368)(16,0.63337746)(17,0.63183450)(18,0.63043918)(19,0.62910477)(20,0.62818513)
        };
        \addplot[forget plot,color=red,mark=square,solid] coordinates {
          (0,0.77163194)(1,0.74343239)(2,0.72306373)(3,0.70711609)(4,0.69515869)(5,0.68507506)(6,0.67740215)(7,0.66974213)(8,0.66339374)(9,0.65866856)(10,0.65423163)(11,0.65069620)(12,0.64773056)(13,0.64568481)(14,0.64382958)(15,0.64212485)(16,0.64073531)(17,0.63919358)(18,0.63771667)(19,0.63688802)(20,0.63582151)
        };
        \addplot[forget plot,color=green,mark=o,solid] coordinates {
          (0,0.63103603)(1,0.63103603)(2,0.63103603)(3,0.63103603)(4,0.63103603)(5,0.63103603)(6,0.63103603)(7,0.63103603)(8,0.63103603)(9,0.63103603)(10,0.63103603)(11,0.63103603)(12,0.63103603)(13,0.63103603)(14,0.63103603)(15,0.63103603)(16,0.63103603)(17,0.63103603)(18,0.63103603)(19,0.63103603)(20,0.63103603)
        };
        \addplot[forget plot,color=blue,mark=x,dashed,mark options={solid}] coordinates {
          (0,0.80182783)(1,0.76977223)(2,0.74555447)(3,0.72779786)(4,0.71320362)(5,0.70160700)(6,0.69168839)(7,0.68378100)(8,0.67671207)(9,0.67052789)(10,0.66606579)(11,0.66159930)(12,0.65770436)(13,0.65451470)(14,0.65158939)(15,0.64965836)(16,0.64746183)(17,0.64558406)(18,0.64389951)(19,0.64255969)(20,0.64097316)
        };
        \addplot[forget plot,color=red,mark=+,dashed,mark options={solid}] coordinates {
          (0,0.80182783)(1,0.76944804)(2,0.74636605)(3,0.72857754)(4,0.71420242)(5,0.70267448)(6,0.69361606)(7,0.68510353)(8,0.67868677)(9,0.67373276)(10,0.66993760)(11,0.66664853)(12,0.66384755)(13,0.66085539)(14,0.65882142)(15,0.65722412)(16,0.65491975)(17,0.65383108)(18,0.65305331)(19,0.65313048)(20,0.65240507)
        };
        \addlegendimage{blue,mark=triangle,solid}
        \addlegendentry{CSL-tri, $m=10$}
        \addlegendimage{blue,mark=x,dashed,mark options={solid}}
        \addlegendentry{CSL-tri, $m=20$}
        \addlegendimage{red,mark=square,solid}
        \addlegendentry{CSL-med, $m=10$}
        \addlegendimage{red,mark=+,dashed,mark options={solid}}
        \addlegendentry{CSL-med, $m=20$}
        \addlegendimage{green,mark=o,solid}
        \addlegendentry{Centralize}
      \end{axis}
    \end{tikzpicture}}
  \hspace{0.2cm}%
  \subfigure{
    \begin{tikzpicture}
      \tikzstyle{every node}=[font=\tiny]
      \pgfplotsset{every axis legend/.append style={at={(0.5,1.20)},anchor=south}, every axis y label/.append style={at={(-0.08,0.5)}}}
      \begin{axis}[title={$\alpha=0.3$}, xtick={0,4,8,12,16,20}, xlabel={$T$}, ylabel={Logloss}, ymin=0.6, ymax=0.84, legend columns=1, legend style={font=\footnotesize}, legend pos=north east, font=\small, height=5.5cm, width=7.9cm, mark repeat=4]
        \addplot[forget plot,color=blue,mark=triangle,solid] coordinates {
          (0,0.77163194)(1,0.74588326)(2,0.72697344)(3,0.71194397)(4,0.70059956)(5,0.69005869)(6,0.68220385)(7,0.67536934)(8,0.66902259)(9,0.66350549)(10,0.65899083)(11,0.65539959)(12,0.65148558)(13,0.64876999)(14,0.64591492)(15,0.64367453)(16,0.64175852)(17,0.64028015)(18,0.63856645)(19,0.63750124)(20,0.63682042)
        };
        \addplot[forget plot,color=red,mark=square,solid] coordinates {
          (0,0.77163194)(1,0.74607459)(2,0.72684646)(3,0.71185944)(4,0.70038219)(5,0.69078232)(6,0.68310220)(7,0.67631042)(8,0.67079283)(9,0.66570869)(10,0.66151479)(11,0.65790039)(12,0.65453663)(13,0.65181877)(14,0.64965490)(15,0.64832308)(16,0.64686844)(17,0.64589042)(18,0.64452852)(19,0.64407303)(20,0.64355451)
        };
        \addplot[forget plot,color=green,mark=o,solid] coordinates {
          (0,0.63103603)(1,0.63103603)(2,0.63103603)(3,0.63103603)(4,0.63103603)(5,0.63103603)(6,0.63103603)(7,0.63103603)(8,0.63103603)(9,0.63103603)(10,0.63103603)(11,0.63103603)(12,0.63103603)(13,0.63103603)(14,0.63103603)(15,0.63103603)(16,0.63103603)(17,0.63103603)(18,0.63103603)(19,0.63103603)(20,0.63103603)
        };
        \addplot[forget plot,color=blue,mark=x,dashed,mark options={solid}] coordinates {
          (0,0.80182783)(1,0.77360990)(2,0.75153301)(3,0.73452166)(4,0.72163694)(5,0.71105430)(6,0.70143209)(7,0.69344760)(8,0.68679302)(9,0.68115060)(10,0.67622818)(11,0.67209382)(12,0.66906458)(13,0.66599776)(14,0.66395471)(15,0.66178032)(16,0.65925460)(17,0.65788382)(18,0.65667987)(19,0.65455374)(20,0.65303674)
        };
        \addplot[forget plot,color=red,mark=+,dashed,mark options={solid}] coordinates {
          (0,0.80182783)(1,0.77379771)(2,0.75255893)(3,0.73527221)(4,0.72283054)(5,0.71212279)(6,0.70422413)(7,0.69807159)(8,0.69213971)(9,0.68722305)(10,0.68362937)(11,0.68085098)(12,0.67825792)(13,0.67607437)(14,0.67368127)(15,0.67245417)(16,0.66919456)(17,0.66828321)(18,0.66656814)(19,0.66692666)(20,0.66587124)
        };
        \addlegendimage{blue,mark=triangle,solid}
        \addlegendentry{CSL-tri, $m=10$}
        \addlegendimage{blue,mark=x,dashed,mark options={solid}}
        \addlegendentry{CSL-tri, $m=20$}
        \addlegendimage{red,mark=square,solid}
        \addlegendentry{CSL-med, $m=10$}
        \addlegendimage{red,mark=+,dashed,mark options={solid}}
        \addlegendentry{CSL-med, $m=20$}
        \addlegendimage{green,mark=o,solid}
        \addlegendentry{Centralize}
      \end{axis}
    \end{tikzpicture}}
  \caption{Test Logloss trajectories on Connect-4. The ordinary mean is omitted from the trajectory panels because it diverges in the main $m=10,\alpha=0.2$ setting.}
  \label{fig:real-connect4-logloss-step4}
\end{figure}
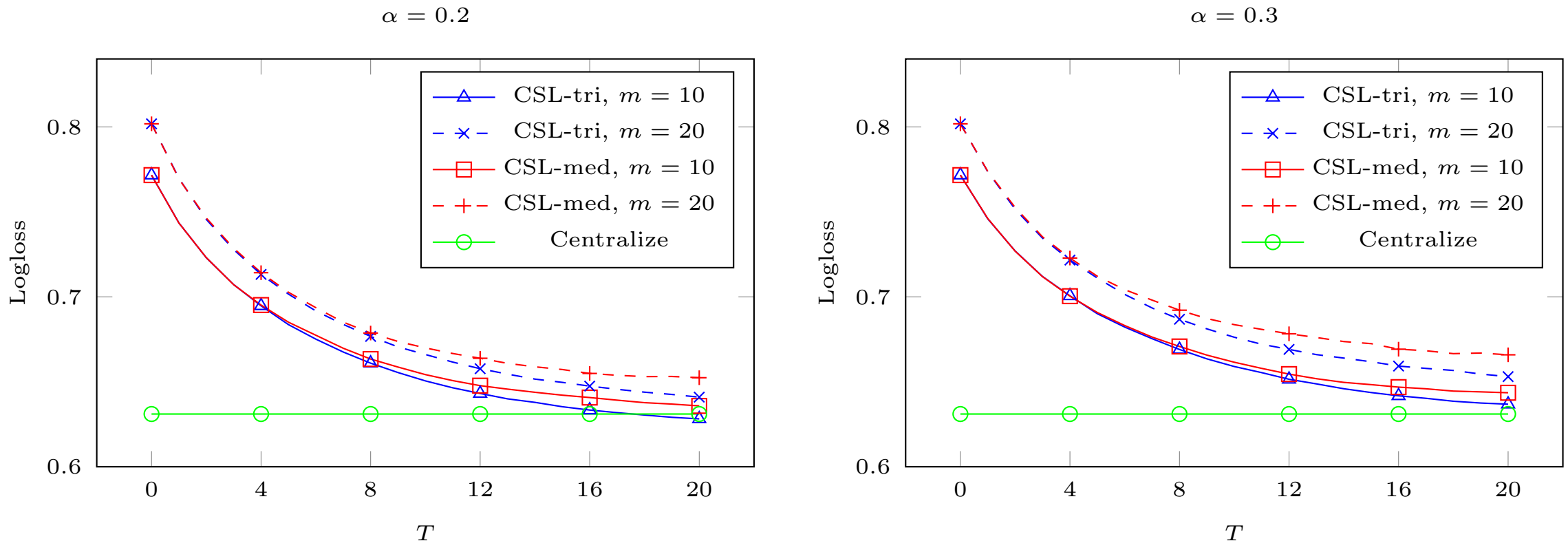

\begin{table}[!ht]
\centering
\caption{Final test performance on Connect-4 for $m=10$ and $\alpha=0.2$. Values are mean $\pm$ standard deviation over 20 repeat  means; each repeat averages five folds.}
\label{tab:real-connect4-main}
\resizebox{\linewidth}{!}{%
\begin{tabular}{lccc}
\toprule
Method & Accuracy & Logloss & Runtime (s) \\
\midrule
Local initialization ($t=0$) & $0.6677\pm0.0020$ & $0.7716\pm0.0023$ & $0.31\pm0.06$ \\
Mean & $0.2677\pm0.1158$ & $20.2333\pm3.1985$ & $15.41\pm0.18$ \\
CSL-tri  & $0.7512\pm0.0017$ & $0.6282\pm0.0035$ & $15.33\pm0.15$ \\
CSL-med  & $0.7507\pm0.0024$ & $0.6358\pm0.0039$ & $15.33\pm0.16$ \\
Centralize & $0.7449\pm0.0002$ & $0.6310\pm0.0001$ & $6.88\pm0.16$ \\
\bottomrule
\end{tabular}%
}
\end{table}

The results separate the quality of the local initialization, the pooled-data reference, and the robustness of distributed aggregation. In the main setting, both robust aggregators improve over the local initialization and remain close to the centralized reference, whereas the ordinary mean is destabilized by the malicious reports. This supports the theoretical role of bounded-influence aggregation: robust methods control the contamination-induced update error, while the ordinary mean permits corrupted gradients to dominate.

Across the trajectory panels, the robust methods improve accuracy and reduce logloss over the outer iterations before approaching a stable regime. Increasing either the Byzantine fraction or the number of workers causes a moderate degradation, but does not produce the catastrophic behavior observed for the ordinary mean. Trimmed mean and median show the same qualitative robustness, with small differences reflecting the coordinate-wise gradient distribution. These observations support the predicted stability ordering, while remaining empirical evidence for one data set and one attack construction rather than a replacement for the theoretical guarantee.

\section{Discussion}\label{SecDiscussion}

The error bound for the robust estimator separates sampling error, governed by the total sample size $N$, from the additional error caused by Byzantine reports, governed by the contamination fraction and local sample size $n$. This separation has a direct implication: increasing the number of workers at fixed $N$ can make robustness harder, even though more workers contribute data. When the Byzantine level $\alpha\asymp m^{-1/2}$, the contamination term has the same dependence on $N$, $s$, and $d$ as the centralized term for fixed $K$; robustness preserves the statistical rate in this regime without eliminating the cost of contamination. Accordingly, the simulations show that the initial improvement gives way to a residual error that reflects this contamination cost. Consistent with this, on real data robust aggregation maintains predictive performance under the tested attack whereas mean aggregation declines.

Together, these guarantees tie two strands of prior work to the sparse multinomial benchmark of \citet{tian}: surrogate-likelihood methods reduce communication through gradient exchange \citep{Communication3}, while robust aggregation controls corrupted gradients \citep{stat}. Our bound inherits both strands, depending on $s\log(Kd)$ rather than on the ambient dimension that appears in a general guarantee \citep{zhou2023}. On top of this, the support-recovery results are a distinct contribution: a small estimation error by itself does not rule out false selected coefficients; Theorem~\ref{theorem2} requires incoherence of the full softmax Hessian to do so, while exact support and sign recovery further require a minimum signal (Theorem~\ref{theorem2prime}). This multiclass condition accounts for cross-class curvature induced by the shared softmax normalization. It also distinguishes the present entrywise support result from the group-penalized discriminant model of \citet{Bao}, whose Gaussian shared-covariance structure does not yield the same condition for conditional softmax regression.

The guarantees require independent sub-Gaussian designs, sufficiently large local samples, and a trimming level that covers the Byzantine fraction. Accordingly, several directions remain open: adapting aggregation to unknown $\alpha$, extending the analysis to heterogeneous workers and broader attacks, and improving the dependence on the number of classes.

\appendix
\section{\textbf{Appendix}}

\subsection{Proof of Theorem~\ref{theorem1}}\label{appA1}
\begin{proof}
	Denote $\widehat{\mathcal{L}}\left(\boldsymbol{\Theta}\right):=\widehat{L}_{1}(\boldsymbol{\Theta})-\left\langle \nabla \widehat{L}_{1}(\widehat{\boldsymbol{\Theta}}_{1})-\nabla \widetilde{\mathcal{L}}(\widehat{\boldsymbol{\Theta}}_{0}),\boldsymbol{\Theta}\right\rangle$ and $\widehat{\boldsymbol{\Delta}}:=\widehat{\boldsymbol{\Theta}}_{1}-\boldsymbol{\Theta}^{\star}$.
	Because $\widehat{\boldsymbol{\Theta}}_{1}$ is the minimum point in~\eqref{eq3}, we have
	\begin{equation}\label{eq5}
		\widehat{\mathcal{L}}\left(\widehat{\boldsymbol{\Theta}}_{1}\right)-\widehat{\mathcal{L}}\left(\boldsymbol{\Theta}^{\star}\right)
		+\lambda\left\|\widehat{\boldsymbol{\Theta}}_{1}\right\|_{1,1}
		-\lambda\left\|\boldsymbol{\Theta}^{\star}\right\|_{1,1} \leq 0.
	\end{equation}
	Since
	$
		\widehat{L}_{1}(\widehat{\boldsymbol{\Theta}}_{1})-\widehat{L}_{1}(\boldsymbol{\Theta}^{\star})- \left\langle \nabla \widehat{L}_{1}(\boldsymbol{\Theta}^{\star}),\widehat{\boldsymbol{\Delta}}\right\rangle\geq0,
	$
	by choosing $\lambda_{1}= 2\|\nabla \widehat{\mathcal{L}}(\boldsymbol{\Theta}^{\star})\|_{\max}$, we obtain that
	\begin{equation}\label{eq6}
	\begin{aligned}
		\widehat{\mathcal{L}}\left(\widehat{\boldsymbol{\Theta}}_{1}\right)-\widehat{\mathcal{L}}\left(\boldsymbol{\Theta}^{\star}\right)&=\widehat{L}_{1}(\widehat{\boldsymbol{\Theta}}_{1})-\widehat{L}_{1}(\boldsymbol{\Theta}^{\star})-
		\left\langle \nabla \widehat{L}_{1}(\widehat{\boldsymbol{\Theta}}_{0})-\nabla \widetilde{\mathcal{L}}(\widehat{\boldsymbol{\Theta}}_{0}),\widehat{\boldsymbol{\Delta}}\right\rangle\\& \geq \left\langle \nabla \widehat{L}_{1}(\boldsymbol{\Theta}^{\star}),\widehat{\boldsymbol{\Delta}}\right\rangle-\left\langle \nabla \widehat{L}_{1}(\widehat{\boldsymbol{\Theta}}_{0})-\nabla \widetilde{\mathcal{L}}(\widehat{\boldsymbol{\Theta}}_{0}),\widehat{\boldsymbol{\Delta}}\right\rangle =\left\langle \nabla \widehat{\mathcal{L}}(\boldsymbol{\Theta}^{\star}),\widehat{\boldsymbol{\Delta}}\right\rangle\\& \geq -\left\| \nabla \widehat{\mathcal{L}}(\boldsymbol{\Theta}^{\star})\right\|_{\text{max}}\left\|\widehat{\boldsymbol{\Delta}}\right\|_{1,1}=-\frac{\lambda}{2}\left\|\widehat{\boldsymbol{\Delta}}\right\|_{1,1}.
	\end{aligned}
    \end{equation}
	Therefore, by combining~\eqref{eq5} with~\eqref{eq6}, we have
	$-\frac{\lambda}{2}\left\|\widehat{\boldsymbol{\Delta}}\right\|_{1,1}\leq 
	\lambda\left\|\boldsymbol{\Theta}^{\star}\right\|_{1,1}-\lambda\left\|\widehat{\boldsymbol{\Theta}}_{1}\right\|_{1,1}$.\par
	Since $\left\|\widehat{\boldsymbol{\Delta}}_{S^{c}}+\boldsymbol{\Theta}^{\star}\right\|_{1,1}=\left\|\widehat{\boldsymbol{\Delta}}_{S^{c}}\right\|_{1,1}+\left\|\boldsymbol{\Theta}^{\star}\right\|_{1,1}$, we have 
	$$
	\begin{aligned}
		-\frac{1}{2}\left\|\widehat{\boldsymbol{\Delta}}_{S}\right\|_{1,1}-\frac{1}{2}\left\|\widehat{\boldsymbol{\Delta}}_{S^{c}}\right\|_{1,1}&\leq 
		\left\|\boldsymbol{\Theta}^{\star}\right\|_{1,1}-\left\|\widehat{\boldsymbol{\Delta}}_{S}+\widehat{\boldsymbol{\Delta}}_{S^{c}}+\boldsymbol{\Theta}^{\star}\right\|_{1,1}\\&\leq
		\left\|\boldsymbol{\Theta}^{\star}\right\|_{1,1}-\left\|\widehat{\boldsymbol{\Delta}}_{S^{c}}+\boldsymbol{\Theta}^{\star}\right\|_{1,1}+\left\|\widehat{\boldsymbol{\Delta}}_{S}\right\|_{1,1}= \left\|\widehat{\boldsymbol{\Delta}}_{S}\right\|_{1,1}-\left\|\widehat{\boldsymbol{\Delta}}_{S^{c}}\right\|_{1,1}
	\end{aligned}
	$$
	Therefore, $\left\|\widehat{\boldsymbol{\Delta}}_{S^{c}}\right\|_{1,1}\leq 3 \left\|\widehat{\boldsymbol{\Delta}}_{S}\right\|_{1,1}$.
	We define $\overline{\boldsymbol{\Theta}}:=\boldsymbol{\Theta}^{\star}+\ell\widehat{\boldsymbol{\Delta}}$, where
	$$
	\ell=\begin{cases}
		1, &\text{if $\|\widehat{\boldsymbol{\Delta}}\|_{F} \leq 1$,}\\
		1/\|\widehat{\boldsymbol{\Delta}}\|_{F}, &\text{if $\|\widehat{\boldsymbol{\Delta}}\|_{F} > 1$.}
	\end{cases}
	$$
	Hence, $\left\|\overline{\boldsymbol{\Theta}}-\boldsymbol{\Theta}^{\star}\right\|_{F}\leq 1$.
	By Lemma 1 in \citet{tian},
	we obtain that
	\begin{align*}
	\mathbb{P}\bigg(\widehat{L}_{1}(\overline{\boldsymbol{\Theta}})-\widehat{L}_{1}(\boldsymbol{\Theta}^{\star})&-\langle\nabla \widehat{L}_{1}(\boldsymbol{\Theta}^{\star}), \overline{\boldsymbol{\Delta}}\rangle \\
	&\geq \tfrac{c_{1}}{K^2}\|\overline{\boldsymbol{\Delta}}\|_2^2- \tfrac{c_{2}}{K^2}\sqrt{\tfrac{\log(d)}{n}}\left\|\overline{\boldsymbol{\Delta}}\right\|_{F}\left\|\overline{\boldsymbol{\Delta}}\right\|_{1,1}\bigg)\geq 1- c_{3}(Kp)^{-1},
	\end{align*}
	where $\overline{\boldsymbol{\Delta}}:=\overline{\boldsymbol{\Theta}}-\boldsymbol{\Theta}^{\star}$. Since $\widehat{\boldsymbol{\Theta}}\in \mathcal{C}(S,3)$, we have  $\overline{\boldsymbol{\Theta}}\in \mathcal{C}(S,3)$. Therefore, as long as $n\gtrsim Ks\log(d)$,
	$\widehat{L}_{1}(\overline{\boldsymbol{\Theta}})-\widehat{L}_{1}(\boldsymbol{\Theta}^{\star})-\langle\nabla \widehat{L}_{1}(\boldsymbol{\Theta}^{\star}), \overline{\boldsymbol{\Delta}}\rangle\geq \frac{c}{K^2}\|\overline{\boldsymbol{\Delta}}\|_{F}^2$.
    Thus,
	$$
	\begin{aligned}
		\widehat{\mathcal{L}}\left(\overline{\boldsymbol{\Theta}}\right)-\widehat{\mathcal{L}}\left(\boldsymbol{\Theta}^{\star}\right)&=\widehat{L}_{1}(\overline{\boldsymbol{\Theta}})-\widehat{L}_{1}(\boldsymbol{\Theta}^{\star})-
		\left\langle \nabla \widehat{L}_{1}(\widehat{\boldsymbol{\Theta}}_{0})-\nabla \widetilde{\mathcal{L}}\left(\widehat{\boldsymbol{\Theta}}_{0}\right),\overline{\boldsymbol{\Delta}}\right\rangle\\& \geq \left\langle \nabla \widehat{\mathcal{L}}(\boldsymbol{\Theta}^{\star}),\overline{\boldsymbol{\Delta}}\right\rangle +\kappa_{\ell}\|\overline{\boldsymbol{\Delta}}\|_{F}^2.
	\end{aligned}
	$$
    Since
	$
	\widehat{\mathcal{L}}\left(\overline{\boldsymbol{\Theta}}\right)-\widehat{\mathcal{L}}\left(\boldsymbol{\Theta}^{\star}\right)
	+\lambda\left\|\overline{\boldsymbol{\Theta}}\right\|_{1,1}
	-\lambda\left\|\boldsymbol{\Theta}^{\star}\right\|_{1,1} \leq 0,
	$
    we obtain that
	$$
	\begin{aligned}
	\lambda\left\|\boldsymbol{\Theta}^{\star}\right\|_{1,1}
	-\lambda\left\|\overline{\boldsymbol{\Theta}}\right\|_{1,1} &\geq \left\langle \nabla \widehat{\mathcal{L}}(\boldsymbol{\Theta}^{\star}),\overline{\boldsymbol{\Delta}}\right\rangle +\kappa_{\ell}\|\overline{\boldsymbol{\Delta}}\|_{F}^2\\&\geq -\|\nabla \widehat{\mathcal{L}}(\boldsymbol{\Theta}^{\star})\|_{\max}\|\overline{\boldsymbol{\Delta}}\|_{1,1}+\kappa_{\ell}\|\overline{\boldsymbol{\Delta}}\|_{F}^2\\&= -\frac{\lambda_{1}}{2}\|\overline{\boldsymbol{\Delta}}\|_{1,1}+\kappa_{\ell}\|\overline{\boldsymbol{\Delta}}\|_{F}^2.
	\end{aligned}
	$$
	Because $	\|\overline{\boldsymbol{\Delta}}\|_{1,1}\geq\left\|\boldsymbol{\Theta}^{\star}\right\|_{1,1}
	-\left\|\overline{\boldsymbol{\Theta}}\right\|_{1,1}$ and $\overline{\boldsymbol{\Theta}}\in \mathcal{C}(S,3)$,
	$$
	\kappa_{\ell}\left\|\overline{\boldsymbol{\Delta}}\right\|_{F}^2\leq \frac{3\lambda_{1}}{2}\left\|\overline{\boldsymbol{\Delta}}\right\|_{1,1}\leq 6\lambda_{1}\left\|\overline{\boldsymbol{\Delta}}_{S}\right\|_{1,1}\leq 6\sqrt{(K-1)s}\lambda_{1}\left\|\overline{\boldsymbol{\Delta}}\right\|_{F}.
	$$
	Therefore, we obtain that $\left\|\overline{\boldsymbol{\Delta}}\right\|_{F}\lesssim \sqrt{Ks}\lambda_{1}$. By Lemma~\ref{lemma2}, as long as $\log(d)^3K^{17/2}s^{5/2}\log (Kd)\lesssim n$,
	$$
	\begin{aligned}
	\left\|\overline{\boldsymbol{\Delta}}\right\|_{F}\lesssim \sqrt{(K-1)s}\lambda_{1}&\lesssim K^{11/2}s^{3/2}\frac{\log (Kd)}{n}+\log(d)^3K^{17/2}s^{5/2}\frac{\log (Kd)}{n}\\&+\frac{\sqrt{1-\alpha}}{1-2\beta}\sqrt{\frac{\log(Kd)}{N}}+\frac{2\beta+\alpha}{1-2\beta}\sqrt{\frac{\log(Kd)}{n}}<1.
	\end{aligned}
	$$ 
	Note that as long as $\left\|\overline{\boldsymbol{\Delta}}\right\|_{F}<1$, we have $\ell=1$. Therefore, we have $\overline{\boldsymbol{\Theta}}=\widehat{\boldsymbol{\Theta}}_{1}$ and 
	$$
	\mathbb{P}\left(\left\|\widehat{\boldsymbol{\Theta}}_{1}-\boldsymbol{\Theta}^{\star}\right\|_{F}\lesssim \frac{\sqrt{1-\alpha}}{1-2\beta}\sqrt{\frac{\log(Kd)}{N}}+\frac{2\beta+\alpha}{1-2\beta}\sqrt{\frac{\log(Kd)}{n}}\right)\geq 1-cmd^{-C}.
	$$
	By math induction, for $\forall t\geq 1$,
    $$
    \begin{aligned}
    \left\|\widehat{\boldsymbol{\Theta}}_{t}-\boldsymbol{\Theta}^{\star}\right\|_{F}&\lesssim \ldots\lesssim\left\|\widehat{\boldsymbol{\Theta}}_{2}-\boldsymbol{\Theta}^{\star}\right\|_{F}\lesssim_{\mathbb{P}} \frac{\sqrt{1-\alpha}}{1-2\beta}\sqrt{\frac{\log(Kd)}{N}}+\frac{2\beta+\alpha}{1-2\beta}\sqrt{\frac{\log(Kd)}{n}}\\ &+\sqrt{\frac{\log (d)}{n}}\left\|\boldsymbol{\Theta}^{\star}-\widehat{\boldsymbol{\Theta}}_{1}\right\|_{1,1}+\log(d)^3 \left\|\boldsymbol{\Theta}^{\star}-\widehat{\boldsymbol{\Theta}}_{1}\right\|_{1,1}^2\\ &\lesssim \frac{\sqrt{1-\alpha}}{1-2\beta}\sqrt{\frac{\log(Kd)}{N}}+\frac{2\beta+\alpha}{1-2\beta}\sqrt{\frac{\log(Kd)}{n}}.
    \end{aligned}
    $$
\end{proof}

\subsection{Technical lemmas}\label{appA2}

\begin{lemma}[\texorpdfstring{$\ell_\infty$ bound on the surrogate gradient at $\boldsymbol{\Theta}^{\star}$}{l-infty bound on the surrogate gradient at Theta-star}]\label{lemma2}
Under Assumption~\ref{assumption1}, for the surrogate objective $\widehat{\mathcal{L}}$ defined in Algorithm~\ref{Al2}, there exist constants $c_{1},c_{2}>0$ independent of $n,m,d,K$ such that
$$
\bigl\|\nabla\widehat{\mathcal{L}}(\boldsymbol{\Theta}^{\star})\bigr\|_{\max}\;\lesssim\;\rho_{N}
$$
with probability at least $1-c_{1}(Kd)^{-c_{2}}$.
\end{lemma}

\begin{lemma}[Uniform empirical-process bound on a sparse index set]\label{lemma3}
Under Assumption~\ref{assumption1}, If $K^{7}s^{2}\log(Kd)\,d\lesssim n$ and $\mathbb{P}\bigl(\|\widehat{\boldsymbol{\Theta}}_{0}\|_{0}\leq C Ks\bigr)\;\geq\;1-c_{1}(Kd)^{-c_{2}}$, then there exist constants $C_{4},c>0$ independent of $n,m,d,K$ such that 
$$
\begin{aligned}
\frac{1}{(1-2\beta)N}&\left\| \sum_{i\in [n],\ell \in \mathcal{U}}\boldsymbol{x}_i^{(\ell)}\left(\boldsymbol{P}_{i}^{\star(\ell)}-\widehat{\boldsymbol{P}}_{i}^{(\ell)}\right)^{\top} -\mathbb{E}\left[\sum_{i\in [n],\ell \in \mathcal{U}}\boldsymbol{x}_i^{(\ell)}\left(\boldsymbol{P}_{i}^{\star(\ell)}-\widehat{\boldsymbol{P}}_{i}^{(\ell)}\right)^{\top} \right]\right\|_{\max}\\&\leq C\left(\sqrt{\frac{Ks \log (d)}{(1-2\beta)N}}+\frac{Ks \sqrt{\log (d)} \log (d)}{(1-2\beta)N}\right)
\end{aligned}
$$
with probability at least $1-d^{-C}$.
\end{lemma}

\subsubsection*{Proof of Lemma~\ref{lemma2}}
\begin{proof}
	For $\forall \boldsymbol{\Theta}=:(\boldsymbol{\theta}_{1},\ldots, \boldsymbol{\theta}_{K-1}) \in \mathbb{R}^{d\times (K-1)}$, denote $p_j(\boldsymbol{x}_i,\boldsymbol{\Theta}):=\frac{\exp\left(\langle \boldsymbol{\theta}_{j},\boldsymbol{x}_i\rangle\right)}{1+\sum\limits_{i=1}^{K-1}\exp\left(\langle \boldsymbol{\theta}_{i},\boldsymbol{x}_i\rangle\right)}$ and $p_{j}^{\prime}\left(\left\langle\boldsymbol{\theta}_{k}, \boldsymbol{x}_i\right\rangle\right)$ is the derivative of $p_j(\boldsymbol{x}_i,\boldsymbol{\Theta}^{\star})$ with respect to $\left\langle\boldsymbol{\theta}_{k}, \boldsymbol{x}_i\right\rangle$.
	Because
	\[
	\resizebox{\linewidth}{!}{$
	\begin{aligned}
		\left\| \nabla \widehat{\mathcal{L}}(\boldsymbol{\Theta}^{\star})\right\|_{\max}=&\left\|\nabla \widehat{L}_{1}(\boldsymbol{\Theta}^{\star})-\nabla \widehat{L}_{1}\left(\widehat{\boldsymbol{\Theta}}_{0}\right)+\nabla \widetilde{\mathcal{L}}\left(\widehat{\boldsymbol{\Theta}}_{0}\right)\right\|_{\max}\\
		=&\left\|\frac{1}{n} \sum_{i=1}^n\boldsymbol{x}_i \left(\boldsymbol{P}_{i}^{\star}-\widehat{\boldsymbol{P}}_{i}\right)^{\top} +\mathop{\text{Aggregate}}\left\{\widehat{\boldsymbol{G}}^{(\ell)}\left(\widehat{\boldsymbol{\Theta}}_{0}\right): \ell \in [m] \right\}\right\|_{\max}\\\leq&\left\|\frac{1}{n} \sum_{i=1}^{n}\boldsymbol{x}_i\left(\boldsymbol{P}_{i}^{\star}-\widehat{\boldsymbol{P}}_{i}\right)^{\top} -\frac{1}{(1-2\beta)N} \sum_{i\in [n],\ell \in \mathcal{U}}\boldsymbol{x}_i^{(\ell)}\left(\boldsymbol{P}_{i}^{\star(\ell)}-\widehat{\boldsymbol{P}}_{i}^{(\ell)}\right)^{\top} \right\|_{\max}\\+&\frac{1}{(1-2\beta)N}\left\| \sum_{i\in [n],\ell \in \mathcal{U}}\boldsymbol{x}_i^{(\ell)}\left(\boldsymbol{P}_{i}^{\star(\ell)}-\widehat{\boldsymbol{P}}_{i}^{(\ell)}\right)^{\top} -\mathbb{E}\left[\sum_{i\in [n],\ell \in \mathcal{U}}\boldsymbol{x}_i^{(\ell)}\left(\boldsymbol{P}_{i}^{\star(\ell)}-\widehat{\boldsymbol{P}}_{i}^{(\ell)}\right)^{\top} \right]\right\|_{\max}\\+&\left\|\mathop{\text{Aggregate}}\limits_{\ell \in [m]}\left\{\widehat{\boldsymbol{G}}^{(\ell)}\left(\widehat{\boldsymbol{\Theta}}_{0}\right)\right\}-\frac{1}{(1-2\beta)N}\mathbb{E}\left[ \sum_{i\in [n],\ell \in \mathcal{U}}\boldsymbol{x}_i^{(\ell)}\widehat{\boldsymbol{P}}_{i}^{(\ell)\top}\right]\right\|_{\max}:=T_{1}+T_{2}+T_{3}
	\end{aligned}
	$}
	\]
	where $\boldsymbol{P}_{i}^{\star}:=\left(p_1(\boldsymbol{x}_i,\boldsymbol{\Theta}^{\star}),\ldots,p_{K-1}(\boldsymbol{x}_i,\boldsymbol{\Theta}^{\star})\right)^{\top},\widehat{\boldsymbol{P}}_{i}:=\left(p_1(\boldsymbol{x}_i,\widehat{\boldsymbol{\Theta}}_{0}),\ldots,p_{K-1}(\boldsymbol{x}_i,\widehat{\boldsymbol{\Theta}}_{0})\right)^{\top}$.
	By Taylor expansion, we have 
	$$p_j(\boldsymbol{x}_i,\widehat{\boldsymbol{\Theta}}_{0})-p_j(\boldsymbol{x}_i,\boldsymbol{\Theta}^{\star})=\sum_{k=1}^{K-1}p_{j}^{\prime}\left(\left\langle\boldsymbol{\theta}_{k}^{\star}, \boldsymbol{x}_i\right\rangle\right) \left\langle\widehat{\boldsymbol{\theta}}_{k}-\boldsymbol{\theta}_{k}^{\star}, \boldsymbol{x}_i\right\rangle+\sum_{k=1}^{K-1}\frac{p_{j}^{\prime\prime}\left(\widetilde{u}_{k,i}\right)}{2} \left\langle\widehat{\boldsymbol{\theta}}_{k}-\boldsymbol{\theta}_{k}^{\star}, \boldsymbol{x}_i\right\rangle^2$$
	where $\widetilde{u}_{k,i}$ is a number between $\left\langle\widehat{\boldsymbol{\theta}}_{k}, \boldsymbol{x}_i\right\rangle$ and $\left\langle\boldsymbol{\theta}_{k}^{\star}, \boldsymbol{x}_i\right\rangle$. $0<p_{j}^{\prime\prime}\left(\widetilde{u}_{k,i}\right)\leq 1$ and $\widehat{\boldsymbol{\Theta}}_{0}:=\left(\widehat{\boldsymbol{\theta}}_{1},\ldots,\widehat{\boldsymbol{\theta}}_{K-1}\right)$.
	$$
	\resizebox{\linewidth}{!}{$
	\begin{aligned}
		T_{1}&=\max_{j\in[K-1]}\left\|\frac{1}{n} \sum_{i=1}^n\left(p_j(\boldsymbol{x}_i,\boldsymbol{\Theta}^{\star})-p_j(\boldsymbol{x}_i,\widehat{\boldsymbol{\Theta}}_{0})\right) \boldsymbol{x}_i-\frac{1}{(1-2\beta)N} \sum_{i=1}^N\left(p_j(\boldsymbol{x}_i,\boldsymbol{\Theta}^{\star})-p_j(\boldsymbol{x}_i,\widehat{\boldsymbol{\Theta}}_{0})\right) \boldsymbol{x}_i\right\|_{\max}\\&\leq\max_{j\in[K\!-\!1]}\left\|\frac{1}{n} \sum_{i=1}^n\sum_{k=1}^{K-1}p_{j}^{\prime}\left(\left\langle\boldsymbol{\theta}_{k}^{\star}, \boldsymbol{x}_i\right\rangle\right) \left\langle\boldsymbol{\theta}_{k}^{\star}\!-\!\widehat{\boldsymbol{\theta}}_{k}, \boldsymbol{x}_i\right\rangle \boldsymbol{x}_i\!-\!\frac{1}{(1-2\beta)N} \sum_{i=1}^{(1-2\beta)N}\sum_{k=1}^{K-1}p_{j}^{\prime}\left(\left\langle\boldsymbol{\theta}_{k}^{\star}, \boldsymbol{x}_i\right\rangle\right) \left\langle\boldsymbol{\theta}_{k}^{\star}\!-\!\widehat{\boldsymbol{\theta}}_{k}, \boldsymbol{x}_i\right\rangle \boldsymbol{x}_i\right\|_{\max}\\&+\max_{j\in[K-1]}\left\|\frac{1}{n} \sum_{i=1}^n\sum_{k=1}^{K-1}\frac{p_{j}^{\prime\prime}\left(\widetilde{u}_{k,i}\right)}{2} \left\langle\widehat{\boldsymbol{\theta}}_{k}-\boldsymbol{\theta}_{k}^{\star}, \boldsymbol{x}_i\right\rangle^2\boldsymbol{x}_i-\frac{1}{(1-2\beta)N}\sum_{i=1}^{(1-2\beta)N}\sum_{k=1}^{K-1}\frac{p_{j}^{\prime\prime}\left(\widetilde{u}_{k,i}\right)}{2} \left\langle\widehat{\boldsymbol{\theta}}_{k}-\boldsymbol{\theta}_{k}^{\star}, \boldsymbol{x}_i\right\rangle^2\boldsymbol{x}_i\right\|_{\max}\\&\leq\sum_{k=1}^{K-1}\max_{j\in[K-1]}\left\|\frac{1}{n} \sum_{i=1}^np_{j}^{\prime}\left(\left\langle\boldsymbol{\theta}_{k}^{\star}, \boldsymbol{x}_i\right\rangle\right) \left\langle\boldsymbol{\theta}_{k}^{\star}-\widehat{\boldsymbol{\theta}}_{k}, \boldsymbol{x}_i\right\rangle \boldsymbol{x}_i-\frac{1}{(1-2\beta)N}\sum_{i=1}^{(1-2\beta)N}p_{j}^{\prime}\left(\left\langle\boldsymbol{\theta}_{k}^{\star}, \boldsymbol{x}_i\right\rangle\right) \left\langle\boldsymbol{\theta}_{k}^{\star}-\widehat{\boldsymbol{\theta}}_{k}, \boldsymbol{x}_i\right\rangle \boldsymbol{x}_i\right\|_{\max}\\&+\max_{j\in[K-1]}\left\|\sum_{k=1}^{K-1}\frac{1}{n} \sum_{i=1}^n\frac{p_{j}^{\prime\prime}\left(\widetilde{u}_{k,i}\right)}{2} \left\langle\widehat{\boldsymbol{\theta}}_{k}-\boldsymbol{\theta}_{k}^{\star}, \boldsymbol{x}_i\right\rangle^2\boldsymbol{x}_i-\frac{1}{(1-2\beta)N}\sum_{i=1}^{(1-2\beta)N}\sum_{k=1}^{K-1}\frac{p_{j}^{\prime\prime}\left(\widetilde{u}_{k,i}\right)}{2} \left\langle\widehat{\boldsymbol{\theta}}_{k}-\boldsymbol{\theta}_{k}^{\star}, \boldsymbol{x}_i\right\rangle^2\boldsymbol{x}_i\right\|_{\max}
		\\&\leq \sum_{k=1}^{K-1}\max_{j\in[K-1]}\left\|\frac{1}{n} \sum_{i=1}^np_{j}^{\prime}\left(\left\langle\boldsymbol{\theta}_{k}^{\star}, \boldsymbol{x}_i\right\rangle\right) \boldsymbol{x}_{i} \boldsymbol{x}_{i}^{\top} -\mathbb{E}[p_{j}^{\prime}\left(\left\langle\boldsymbol{\theta}_{k}^{\star}, \boldsymbol{x}_i\right\rangle\right) \boldsymbol{x}_{i} \boldsymbol{x}_{i}^{\top} ]\right\|_{\max}\left\|\boldsymbol{\theta}_{k}^{\star}-\widehat{\boldsymbol{\theta}}_{k}\right\|_{1}\\&+\sum_{k=1}^{K-1}\max_{j\in[K-1]}\left\|\frac{1}{(1-2\beta)N}\sum_{i=1}^{(1-2\beta)N} p_{j}^{\prime}\left(\left\langle\boldsymbol{\theta}_{k}^{\star}, \boldsymbol{x}_i\right\rangle\right) \boldsymbol{x}_{i} \boldsymbol{x}_{i}^{\top} -\mathbb{E}[p_{j}^{\prime}\left(\left\langle\boldsymbol{\theta}_{k}^{\star}, \boldsymbol{x}_i\right\rangle\right) \boldsymbol{x}_{i} \boldsymbol{x}_{i}^{\top} ]\right\|_{\max}\left\|\boldsymbol{\theta}_{k}^{\star}-\widehat{\boldsymbol{\theta}}_{k}\right\|_{1}\\&+\max_{j\in[K-1]}\left\{\left\|\widehat{\boldsymbol{\theta}}_{k}-\boldsymbol{\theta}_{k}^{\star}\right\|_{1}\right\}\sum_{j=1}^{K-1}\left\|\widehat{\boldsymbol{\theta}}_{k}-\boldsymbol{\theta}_{k}^{\star}\right\|_{1}\|\boldsymbol{x}_i\|_{\max}^3=:T_{11}+T_{12}+T_{13}
	\end{aligned}
	$}
    $$
	Since $\|\boldsymbol{x}_i\|_{\psi_{2}}\leq \kappa_{0}\leq \infty$ and $|p_{j}^{\prime}\left(\left\langle\boldsymbol{\theta}_{k}^{\star}, \boldsymbol{x}_i\right\rangle\right)|\leq 1$, we have $\|p_{j}^{\prime}\left(\left\langle\boldsymbol{\theta}_{k}^{\star}, \boldsymbol{x}_i\right\rangle\right)x_{ij}x_{ik}\|_{\psi_{1}}\leq \kappa_{0}$. By Proposition 5.16
(Bernstein-type inequality) in Vershynin (2010), we have 
    $\mathbb{P}\left(\left\|\boldsymbol{x}_i\right\|_{\max} \geq t\right) \leq 2 d e^{-\frac{t^2}{2 \sigma^2}}$ and with probability at least $1-d^{-\delta}$,
	$$
	\left\|\frac{1}{n} \sum_{i=1}^np_{j}^{\prime}\left(\left\langle\boldsymbol{\theta}_{k}^{\star}, \boldsymbol{x}_i\right\rangle\right) \boldsymbol{x}_{i} \boldsymbol{x}_{i}^{\top} -\mathbb{E}\left[p_{j}^{\prime}\left(\left\langle\boldsymbol{\theta}_{k}^{\star}, \boldsymbol{x}_i\right\rangle\right)\boldsymbol{x}_{i} \boldsymbol{x}_{i}^{\top}\right]\right\|_{\max} \lesssim \sqrt{\frac{\delta\log (d)}{n}},
	$$
	$$
	\left\|\frac{1}{(1-2\beta)N}\sum_{i=1}^{(1-2\beta)N} p_{j}^{\prime}\left(\left\langle\boldsymbol{\theta}_{k}^{\star}, \boldsymbol{x}_i\right\rangle\right)\boldsymbol{x}_{i} \boldsymbol{x}_{i}^{\top} -\mathbb{E}\left[p_{j}^{\prime}\left(\left\langle\boldsymbol{\theta}_{k}^{\star}, \boldsymbol{x}_i\right\rangle\right) \boldsymbol{x}_{i} \boldsymbol{x}_{i}^{\top} \right]\right\|_{\max}\lesssim \sqrt{\frac{\delta \log (d)}{(1-2\beta)N}}.
	$$
    Therefore, as long as $\log(d)^6K^{12}s^4\log(Kd)m\lesssim n$,
	$$
	T_{1}\lesssim \sqrt{\frac{\log (d)}{n}}\left\|\boldsymbol{\Theta}^{\star}-\widehat{\boldsymbol{\Theta}}_{0}\right\|_{1,1}+\log(d)^3 \left\|\boldsymbol{\Theta}^{\star}-\widehat{\boldsymbol{\Theta}}_{0}\right\|_{1,1}^2\lesssim_{\mathbb{P}}\sqrt{\frac{\log (Kd)}{N}}.
	$$
	Since $\mathbb{P}\left(\left\|\boldsymbol{\Theta}^{\star}-\widehat{\boldsymbol{\Theta}}_{0}\right\|_{1,1} \leq  CK^{3}\sqrt{\frac{s\log(Kd)}{n}}\right)\geq 1-c(Kd)^{-C}$ and Theorem~\ref{theorem2}, according Lemma~\ref{lemma3}, we have
	$
	\mathbb{P}\left(T_{2}\leq \sqrt{\frac{\log (d)}{(1-2\beta)N}}+\frac{Ks \sqrt{\log (d)} \log (d)}{(1-2\beta)N}\right)\geq 1-c(Kd)^{-C}.
	$
    \\
	For $T_{3}$, denote $\widehat{g}^{(\ell)}_{j,k}:=\left[\widehat{\boldsymbol{G}}^{(\ell)}\left(\widehat{\boldsymbol{\Theta}}_{0}\right)\right]_{j,k}, \bar{g}_{j,k}:=\frac{1}{(1-2\beta)m} \sum\limits_{\ell \in \mathcal{U}}g_{j,k}^{(\ell)}$, where $g_{j,k}^{(\ell)}:=\mathbb{E}\left[\frac{1}{n} \sum\limits_{i\in [n]}\boldsymbol{x}_i^{(\ell)}\widehat{\boldsymbol{P}}_{i}^{(\ell)\top}\right]_{j,k}$ and $r_{j,k}:=\operatorname{argmax}_{\ell \in \mathcal{B} \cap \mathcal{U}}\left|\widehat{g}^{(\ell)}_{j,k}-\bar{g}_{j,k}\right|$. On the one hand, if $\widehat{g}^{(r_{j,k})}_{j,k}-\bar{g}_{j,k}\geq 0$,
	since $\beta\geq \alpha$, we have $\exists i_{j,k} \in \mathcal{T}\cap \mathcal{M}$, s.t. $ \widehat{g}^{(i_{j,k})}_{j,k}\geq \widehat{g}^{(r_{j,k})}_{j,k}$. Therefore, $\left|\widehat{g}^{(i_{j,k})}_{j,k}-\bar{g}_{j,k}\right|\geq \left|\widehat{g}^{(r_{j,k})}_{j,k}-\bar{g}_{j,k}\right|$. On the other hand, if $\widehat{g}^{(r_{j,k})}_{j,k}-\bar{g}_{j,k}< 0$, $\exists i \in \mathcal{T}\cap \mathcal{M}$, s.t. $ \widehat{g}^{(i_{j,k})}_{j,k}\leq \widehat{g}^{(r_{j,k})}_{j,k}$. Therefore, $\left|\widehat{g}^{(i_{j,k})}_{j,k}-\bar{g}_{j,k}\right|\geq \left|\widehat{g}^{(r_{j,k})}_{j,k}-\bar{g}_{j,k}\right|$. Thus,
	\begin{equation}\label{eq11}
	\resizebox{\linewidth}{!}{$
	\begin{aligned}
		&\left|\sum_{\ell \in \mathcal{B} \cap \mathcal{U}}\left(\widehat{g}^{(\ell)}_{j,k}-\bar{g}_{j,k}\right)\right|\leq\alpha m \max_{\ell \in \mathcal{B} \cap \mathcal{U}}\left|\widehat{g}^{(\ell)}_{j,k}-\bar{g}_{j,k}\right|\leq \alpha m  \max_{\ell \in \mathcal{T}\cap \mathcal{M}}\left|\widehat{g}^{(\ell)}_{j,k}-\bar{g}_{j,k}\right|\leq \alpha m  \max_{\ell \in  \mathcal{M}}\left|\widehat{g}^{(\ell)}_{j,k}-\bar{g}_{j,k}\right|.
	\end{aligned}
	$}
	\end{equation}
	By~\eqref{eq11} and Lemma 3 of \citet{stat}, we have with probability at least $1-(1-\alpha)me^{-t}-e^{-t}$, 
	$$
	\begin{aligned}
		&\left|\operatorname{trmean}_{\beta}\left\{\widehat{g}^{(\ell)}_{j,k}: \ell \in [m]\right\}-\bar{g}_{j,k}\right|  =\left|\frac{1}{(1-2 \beta) m} \sum_{\ell \in \mathcal{U}} \widehat{g}^{(\ell)}_{j,k}-\bar{g}_{j,k}\right| \\
		=&\frac{1}{(1-2 \beta) m}\left|\sum_{\ell \in \mathcal{M}}\left(\widehat{g}^{(\ell)}_{j,k}-g^{(\ell)}_{j,k}\right)-\sum_{\ell \in \mathcal{M} \cap \mathcal{T}}\left(\widehat{g}^{(\ell)}_{j,k}-g^{(\ell)}_{j,k}\right)+\sum_{\ell \in \mathcal{B} \cap \mathcal{U}}\left(\widehat{g}^{(\ell)}_{j,k}-g^{(\ell)}_{j,k}\right)\right| \\
		\leq&\frac{1}{(1-2 \beta) m}\left(\left|\sum_{\ell \in \mathcal{M}}\left(\widehat{g}^{(\ell)}_{j,k}-g^{(\ell)}_{j,k}\right)\right|+\left|\sum_{\ell \in \mathcal{M} \cap \mathcal{T}}\left(\widehat{g}^{(\ell)}_{j,k}-g^{(\ell)}_{j,k}\right)\right|+\left|\sum_{\ell \in \mathcal{B} \cap \mathcal{U}}\left(\widehat{g}^{(\ell)}_{j,k}-g^{(\ell)}_{j,k}\right)\right|\right)\\
		\leq&\frac{1}{(1-2 \beta) m}\left(\left|\sum_{\ell \in \mathcal{M}}\left(\widehat{g}^{(\ell)}_{j,k}-g^{(\ell)}_{j,k}\right)\right|+2\beta m \max_{\ell \in \mathcal{M}}\left|\widehat{g}^{(\ell)}_{j,k}-g^{(\ell)}_{j,k}\right|+\alpha m  \max_{\ell \in  \mathcal{M}}\left|\widehat{g}^{(\ell)}_{j,k}-g^{(\ell)}_{j,k}\right|\right)\\
		\leq&\frac{1}{(1-2 \beta) m}\left(\left|\sum_{\ell \in \mathcal{M}}\left(\widehat{g}^{(\ell)}_{j,k}-g^{(\ell)}_{j,k}\right)\right|+(2\beta+\alpha) m \max_{\ell \in \mathcal{M}}\left|\widehat{g}^{(\ell)}_{j,k}-g^{(\ell)}_{j,k}\right|\right)=:T_{31}+T_{32}
	\end{aligned}
	$$
	Furthermore,
	$$
	\resizebox{\linewidth}{!}{$
	\begin{aligned}
		T_{31}:=&\left|\frac{1}{(1-\alpha) m} \sum_{\ell \in[m] \backslash \mathcal{B}}\left[\nabla{\widehat{L}}_{\ell}\left(\widehat{\boldsymbol{\Theta}}_{0}\right)\right]_{j,k}-\mathbb{E}\left[\frac{1}{(1-\alpha) m} \sum_{\ell \in[m] \backslash \mathcal{B}}\nabla{\widehat{L}}_{\ell}\left(\widehat{\boldsymbol{\Theta}}_{0}\right)\right]_{j,k}\right| \\ \leq & \left|\frac{1}{(1-\alpha) m} \sum_{\ell \in[m] \backslash \mathcal{B}}\left[\nabla{\widehat{L}}_{\ell}\left(\widehat{\boldsymbol{\Theta}}_{0}\right)-\nabla \widehat{L}_{\ell}(\boldsymbol{\Theta}^{\star})\right]_{j,k}-\frac{1}{(1-\alpha) m} \sum_{\ell \in[m] \backslash \mathcal{B}}\mathbb{E}\left[\nabla{\widehat{L}}_{\ell}\left(\widehat{\boldsymbol{\Theta}}_{0}\right)-\nabla \widehat{L}_{\ell}(\boldsymbol{\Theta}^{\star})\right]_{j,k}\right|\\+&
		\left|\frac{1}{(1-\alpha) m} \sum_{\ell \in[m] \backslash \mathcal{B}}\left[\nabla \widehat{L}_{\ell}(\boldsymbol{\Theta}^{\star})\right]_{j,k}-\mathbb{E}\left[\nabla \widehat{L}_{\ell}(\boldsymbol{\Theta}^{\star})\right]_{j,k}\right|=:T_{311}+T_{312}.
	\end{aligned}
	$}
	$$
    By following the proof of $T_{2}$, we have $T_{311}\lesssim_{\mathbb{P}}\sqrt{\frac{\log(Kd)}{N}}$.
    For $T_{312}$, since $\left|p_k(\boldsymbol{x}_i,\boldsymbol{\Theta}^{\star})\right|\leq 1, $ $\left\|x_{i,j}p_k(\boldsymbol{x}_i,\boldsymbol{\Theta}^{\star})\right\|_{\psi_{2}}\leq \|x_{i,j}\|_{\psi_{2}}\leq \kappa_{0}$, by Hoeffding-type inequality and union bound, we have that
	$$
    \mathbb{P}\left(\left|\frac{1}{(1-\alpha) m} \sum_{\ell \in[m] \backslash \mathcal{B}}\left[\nabla{\widehat{L}}_{\ell}\left(\boldsymbol{\Theta}^{\star}\right)\right]_{j,k}-\mathbb{E}\left[\nabla{\widehat{L}}_{\ell}\left(\boldsymbol{\Theta}^{\star}\right)\right]_{j,k}\right|\leq \sqrt{\frac{t}{(1-\alpha)nm}}\right)>1-e^{-t}.
    $$
For $T_{32}$, by the same way, we obtain that $T_{32}\lesssim_{\mathbb{P}} \frac{2\beta+\alpha}{1-2\beta}\sqrt{\frac{\log(Kd)}{n}}$.
	Therefore, 
	$$\mathbb{P}\left(T_{3}\leq
	\frac{\sqrt{1-\alpha}}{1-2\beta}\sqrt{\frac{\log(Kd)}{N}}+\frac{2\beta+\alpha}{1-2\beta}\sqrt{\frac{\log(Kd)}{n}}\right)\geq 1-c(Kd)^{-C}-c(1-\alpha)m(Kd)^{-C}.$$
	After the above derivation, we get
    $$
    \left\| \nabla \widehat{\mathcal{L}}(\boldsymbol{\Theta}^{\star})\right\|_{\max}\lesssim_{\mathbb{P}} \sqrt{\frac{\log (Kd)}{N}}+\frac{\sqrt{1-\alpha}}{1-2\beta}\sqrt{\frac{\log(Kd)}{N}}+\frac{2\beta+\alpha}{1-2\beta}\sqrt{\frac{\log(Kd)}{n}}.
    $$
\end{proof}

\subsubsection*{Proof of Lemma~\ref{lemma3}}
The argument follows the proof of Lemma~A.1 in \citet{WangLei} and uses the sparsity condition $\max_{k}\|\widehat{\boldsymbol{\theta}}_{k}\|_{0}\leq C_{2}s$ from Assumption~\ref{assumption1}. Set
$$\Omega=\left\{\boldsymbol{\Theta} \in \mathbb{R}^{d\times (K-1)}:\max\limits_{j\in[K-1]}\left\|\boldsymbol{\theta}_{j}\right\|_0 \leq Cs,\left\|\boldsymbol{\Theta}-\boldsymbol{\Theta}^{\star}\right\|_{1,1} \leq C K^{3}s \sqrt{\log(Kd) / n}\right\},$$
$$
\mathcal{G}_{j}:=\left\{\bigl(p_{k}(\boldsymbol{x},\boldsymbol{\Theta}^{\star})-p_{k}(\boldsymbol{x},\boldsymbol{\Theta})\bigr)[\boldsymbol{x}]_{j}I\{\max_j[\boldsymbol{x}]_j \leq c_d\}:\ \boldsymbol{\Theta}\in\Omega\right\},\quad j\in[d].
$$
with squared integrable envelope function $F(\boldsymbol{x})=\left|[\boldsymbol{x}]_{j}\right|$. We decompose $\Omega$ as $\Omega=\cup_{\substack{T \subset [d]\times\ldots\times [d],\\ |T_j|\leq Cs,\forall j\in [K \!-\!1]}} \Omega(T)$ with $\Omega(T)=\{\boldsymbol{\Theta}\in \Omega$ : support of $\boldsymbol{\theta}_{j} \subset T_{j},\forall j \in [K-1]\} $. We also define $$\mathcal{G}_j(T):=\left\{\left(p_k(\boldsymbol{x},\boldsymbol{\Theta}^{\star})-p_k(\boldsymbol{x}_i,\boldsymbol{\Theta})\right) [\boldsymbol{x}]_{j}I\{\max_j[\boldsymbol{x}]_j \leq c_d\}: \boldsymbol{\Theta} \in \Omega(T)\right\}.$$ 
By Lemma 2.6.15, Lemma 2.6.18 (vi) and (viii) in \citet{Van}, for each fixed $T \subset [d]\times\ldots\times [d]$ with $|T| \leq (K-1)Cs, \mathcal{G}_j(T)$ is a VC-subgraph with index bounded by $(K-1)Cs$ and by Theorem 2.6.7 of \citet{Van}, we have

$$
N\left(\epsilon, \mathcal{G}_j(T), L_2\left(P_{(1-2\beta)N}\right)\right) \leq\left(\frac{C\|F\|_{L_2\left(P_{(1-2\beta)N}\right)}}{\epsilon}\right)^{(K-1)Cs}.
$$
Since there are at most $\binom{d}{Cs}^{(K-1)} \leq (e d /Cs)^{(K-1)Cs}$ different such $T$, we have
$$
N\left(\epsilon, \mathcal{G}_j, L_2\left(P_{(1-2\beta)N}\right)\right) \leq\left(\frac{C\|F\|_{L_2\left(P_{(1-2\beta)N}\right)}}{\epsilon}\right)^{(K-1)Cs}\left(\frac{e d}{Cs}\right)^{(K-1)Cs} \leq\left(\frac{Cd \|F\|_{L_2\left(P_{(1-2\beta)N}\right)}}{\epsilon}\right)^{(K-1)Cs}
$$
and thus
$$
N\left(\epsilon, \cup_{j=1}^d \mathcal{G}_j, L_2\left(P_{(1-2\beta)N}\right)\right) \leq d\left(\frac{C d}{\epsilon}\right)^{(K-1)Cs}.
$$
Let $\sigma^2=\sup _{f \in \cup_k \mathcal{G}_j} P f^2$. Then by Theorem 3.12 of \citet{Koltch}, we have

$$
\mathbb{E}\left\|R_{(1-2\beta)N}\right\|_{\cup \mathcal{G}_j} \leq C\left(\sigma \sqrt{\frac{Ks \log (d)}{(1-2\beta)N}}+\frac{Ks c_{n}\log (d)}{(1-2\beta)N}\right)
$$
where $\left\|R_{(1-2\beta)N}\right\|_{\cup_j \mathcal{G}_j}=\sup _{f \in \cup_j \mathcal{G}_j} \frac{1}{(1-2\beta)N} \sum_{i=1}^{(1-2\beta)N} \varepsilon_i f\left(\boldsymbol{x}_i, y_i\right)$ with $\varepsilon_i$ being i.i.d. Rademacher random variables. Using the symmetrization inequality, which states that $\mathbb{E}\left\|P_{(1-2\beta)N}-P\right\|_{\cup_j \mathcal{G}_j} \leq$ $2 \mathbb{E}\left\|R_{(1-2\beta)N}\right\|_{\cup_j \mathcal{G}_{j}}$, where $\left\|P_{(1-2\beta)N}-P\right\|_{\cup_j \mathcal{G}_j}=\sup _{f \in \cup_j \mathcal{G}_j} \frac{1}{(1-2\beta)N} \sum_i f\left(\boldsymbol{x}_i, y_i\right)-\mathbb{E} f(\boldsymbol{x}, y)$, Talagrand's inequality (Theorem 2.6 of \citet{Koltch}) gives
$$
\mathbb{P}\left(\left\|P_{(1-2\beta)N}-P\right\|_{\cup_j \mathcal{G}_j} \geq C\left(\sigma \sqrt{\frac{Ks \log (d)}{(1-2\beta)N}}+\frac{Ks c_{d}\log (d)}{(1-2\beta)N}+\sqrt{\frac{\sigma^2 t}{(1-2\beta)N}}+\frac{\sqrt{\log (d)} t}{(1-2\beta)N}\right)\right) \leq e^{-t}
$$
$$
\begin{aligned}
	\sigma^2=&\mathbb{E}\left[\left(p_k(\boldsymbol{x},\boldsymbol{\Theta}^{\star})-p_k(\boldsymbol{x}_i,\boldsymbol{\Theta})\right)^2\right]\leq \mathbb{E}\left[\left\langle \nabla p_k(\boldsymbol{x}_{i},\boldsymbol{\Theta}^{\star}),\boldsymbol{\Theta}^{\star}-\boldsymbol{\Theta}\right\rangle^2\right]\\ \leq & \mathbb{E}\left[\left\| \nabla p_k(\boldsymbol{x}_{i},\boldsymbol{\Theta}^{\star})\right\|_{\max}^2\left\|\boldsymbol{\Theta}^{\star}-\boldsymbol{\Theta}\right\|_{1,1}^2\right]\leq cK^6s^2 \frac{\log(Kd)}{n}.
\end{aligned}
$$
as long as $K^7s^3\log(Kd)\lesssim n$, with probability at least $1-d^{-C}$
$$
\left\|P_{(1-2\beta)N}-P\right\|_{\cup_j \mathcal{G}_j} \leq C\left(\sqrt{\frac{\log (d)}{(1-2\beta)N}}+\frac{Ks \sqrt{\log (d)} \log (d)}{(1-2\beta)N}\right).
$$

\subsection{\texorpdfstring{Proof of Theorem~\ref{theorem2}}{Proofs of support recovery theorems}}\label{appA3}

Under an additional minimum-signal condition, the conclusion of Theorem~\ref{theorem2} can be strengthened to exact support recovery and sign consistency, since the support-inclusion part is obtained as the first half of the argument for Theorem~\ref{theorem2prime}.
\begin{assumption}\label{assumption3}
The minimum magnitude of the nonzero entries of $\boldsymbol{\Theta}^{\star}$ satisfies
$$
\Theta^{\star}_{\min}:=\min_{(j,k)\in S}|\Theta^{\star}_{j,k}|>2 D_{\infty}\!\left(1+\tfrac{\gamma}{8}\right)\lambda.
$$
\end{assumption}
Assumption~\ref{assumption3} is the multinomial counterpart of the minimum-value condition used in Theorem~1 of \citet{Ravikumar2010}.

\begingroup
\addtocounter{theorem}{-1}%
\renewcommand{\thetheorem}{\arabic{theorem}$'$}%
\renewcommand{\theHtheorem}{\arabic{theorem}prime}%
\begin{theorem}[Support recovery and sign consistency]\label{theorem2prime}
Under the conditions of Theorem~\ref{theorem2} together with Assumption~\ref{assumption3}, on the same high-probability event

$$
\operatorname{supp}(\widehat{\boldsymbol{\Theta}}_{t})=S,\qquad \operatorname{sign}(\widehat{\Theta}_{t,j,k})=\operatorname{sign}(\Theta^{\star}_{j,k})\quad\forall (j,k)\in S.
$$
\end{theorem}
\endgroup

\begin{proof}[Proof of Theorem~\ref{theorem2prime}]
Recall the notation $S$, $\boldsymbol{H}^{\star}$ and $|||\cdot|||_{\infty}$ from Section~\ref{SecMain}. 
Let $\mathcal{E}_1$ denote the event that $\big\|\nabla \widehat{\mathcal{L}}(\boldsymbol{\Theta}^*)\big\|_{\max} \leq \rho_N$, then by Lemma~\ref{lemma2}, $\mathbb{P}(\mathcal{E}_1) \geq c_1 (K d)^{-c_2}$.
We construct a candidate primal--dual pair $(\check{\boldsymbol{\Theta}},\check{\boldsymbol{Z}})$ as follows. \newline
(i) Let $(\check{\boldsymbol{\Theta}}_{S},\check{\boldsymbol{Z}}_{S})$ be a solution and the associated subgradient of the oracle restricted problem
$$
\min_{\boldsymbol{\Theta}:\,\boldsymbol{\Theta}_{S^{c}}=\boldsymbol{0}}\Bigl\{\widehat{\mathcal{L}}(\boldsymbol{\Theta})+\lambda\|\boldsymbol{\Theta}\|_{1,1}\Bigr\},
$$
so that $\check{\boldsymbol{Z}}_{S}\in\partial\|\check{\boldsymbol{\Theta}}_{S}\|_{1,1}$, $\|\check{\boldsymbol{Z}}_{S}\|_{\infty}\leq 1$, and the stationarity condition on $S$ holds:
\begin{equation}\label{eq:KKT}
[\nabla\widehat{\mathcal{L}}(\check{\boldsymbol{\Theta}})]_{S}+\lambda\check{\boldsymbol{Z}}_{S}=\boldsymbol{0}.
\end{equation}
where $\check{\boldsymbol{\Theta}}:=(\check{\boldsymbol{\Theta}}_{S}:\check{\boldsymbol{\Theta}}_{S^c}=0)$
(ii) Define $\check{\boldsymbol{Z}}_{S^{c}}$ such that the full KKT equation $\nabla\widehat{\mathcal{L}}(\check{\boldsymbol{\Theta}})+\lambda\check{\boldsymbol{Z}}=\boldsymbol{0}$ holds on $S^{c}$, namely
$$
\check{\boldsymbol{Z}}_{S^{c}}:=-\tfrac{1}{\lambda}[\nabla\widehat{\mathcal{L}}(\check{\boldsymbol{\Theta}})]_{S^{c}}.
$$

\noindent\textbf{Claim.} 
Once $\|\check{\boldsymbol{Z}}_{S^{c}}\|_{\infty}<1$, then $\check{\boldsymbol{Z}}\in\partial\|\check{\boldsymbol{\Theta}}\|_{1,1}$ and $\check{\boldsymbol{\Theta}}$ satisfies the KKT condition of the unrestricted problem~\eqref{eq3}.
In addition, strict convexity guarantees that $\check{\boldsymbol\Theta}$ is the unique optimal solution, which yields $\widehat{\boldsymbol\Theta}_t = \check{\boldsymbol\Theta}$ and further implies $\operatorname{supp}(\widehat{\boldsymbol\Theta}_t) \subseteq S$.

Consequently, the proof target reduces to verifying the inequality $\|\check{\boldsymbol Z}_{S^c}\|_\infty < 1$. 
The remaining steps 1--4 establish this strict dual feasibility.

\noindent\textbf{Step 1: Support Set Decomposition.} 
Since $\check{\boldsymbol\Theta}_{S^c} = \boldsymbol\Theta^*_{S^c}$, define
$$
\check{\boldsymbol\Delta} := \operatorname{vec}(\check{\boldsymbol\Theta} - \boldsymbol\Theta^*), \quad \check{\boldsymbol\Delta}_{S^c} = \mathbf{0}.
$$
Given that $\nabla^2 \widehat{\mathcal{L}} \equiv \nabla^2 \widehat{\mathcal{L}}_1$ is smooth, the integral mean value theorem for vector-valued mappings applies:
\begin{equation}\label{eq:mvt}
\nabla\widehat{\mathcal{L}}(\check{\boldsymbol\Theta}) - \nabla\widehat{\mathcal{L}}(\boldsymbol\Theta^*) = \bar{\boldsymbol H} \check{\boldsymbol\Delta}, \quad \bar{\boldsymbol H} := \int_{0}^{1} \nabla^2 \widehat{\mathcal{L}}_1(\boldsymbol\Theta^* + t\check{\boldsymbol\Delta}) dt.
\end{equation}
Substitute~\eqref{eq:mvt} and partition over $S$ and $S^c$ gives
$$
\begin{aligned}
&\bar{\boldsymbol H}_{SS} \check{\boldsymbol\Delta}_S + \big[\nabla\widehat{\mathcal{L}}(\boldsymbol\Theta^*)\big]_S + \lambda \check{\boldsymbol Z}_S = \mathbf{0}, \\
&\bar{\boldsymbol H}_{S^c S} \check{\boldsymbol\Delta}_S + \big[\nabla\widehat{\mathcal{L}}(\boldsymbol\Theta^*)\big]_{S^c} + \lambda \check{\boldsymbol Z}_{S^c} = \mathbf{0}.
\end{aligned}
$$
Solve for $\check{\boldsymbol\Delta}_S$ from block $S$ (the invertibility of $\bar{\boldsymbol H}_{SS}$ is verified in Step 2):
\begin{equation}\label{eq:DeltaS}
\check{\boldsymbol\Delta}_S = -\bar{\boldsymbol H}_{SS}^{-1}\Big(\big[\nabla\widehat{\mathcal{L}}(\boldsymbol\Theta^*)\big]_S + \lambda \check{\boldsymbol Z}_S\Big).
\end{equation}
Substitute~\eqref{eq:DeltaS} into block $S^c$ and rearrange terms to obtain
\begin{equation}\label{eq:ZSc}
\check{\boldsymbol Z}_{S^c} = \bar{\boldsymbol H}_{S^c S}\bar{\boldsymbol H}_{SS}^{-1}\check{\boldsymbol Z}_S + \frac{1}{\lambda}\Big(\bar{\boldsymbol H}_{S^c S}\bar{\boldsymbol H}_{SS}^{-1}\big[\nabla\widehat{\mathcal{L}}(\boldsymbol\Theta^*)\big]_S - \big[\nabla\widehat{\mathcal{L}}(\boldsymbol\Theta^*)\big]_{S^c}\Big).
\end{equation}

\noindent\textbf{Step2: The Invertibility of $\bar{\boldsymbol H}_{SS}$.}
Define $\boldsymbol{E} := \bar{\boldsymbol{H}} - \boldsymbol{H}^\star$.
$$
\begin{aligned}
\|\boldsymbol{E}\|_{\max} \leq \sup_{t\in[0,1]} \bigg\| \nabla^2 \widehat{\mathcal{L}}_1\big(\boldsymbol{\Theta}^\star + t\check{\boldsymbol{\Delta}}\big) - \nabla^2 \widehat{\mathcal{L}}_1(\boldsymbol{\Theta}^\star) \bigg\|_{\max} + \bigg\| \frac{1}{n}\sum_{i} \boldsymbol{A}_i(\boldsymbol{\Theta}^\star) \otimes \boldsymbol{x}_i \boldsymbol{x}_i^\top - \boldsymbol{H}^\star \bigg\|_{\max}.
\end{aligned}
$$
The inequality $0 \preceq \boldsymbol{A}_i \preceq \frac{1}{2}\boldsymbol{I}$ implies $|[\boldsymbol{A}_i]_{kk'}| \leq \frac{1}{2}$. Under the sub-Gaussian design specified in Assumption~\ref{a1}, each term $x_{ij}x_{ij'}$ is sub-exponential with $\|x_{ij}x_{ij'}\|_{\psi_1} \leq \nu^2$. We apply the entrywise sub-exponential Bernstein inequality and perform the union bound over all $(Kd)^2$ matrix entries.
Let $\mathcal{E}_2$ denote the event that $
R_n := \max_{i\le n}\|\boldsymbol{x}_i\|_\infty \le c\sqrt{\log(Kd)}, \max_{j,j'}\frac{1}{n}\sum_{i}|x_{ij}x_{ij'}|\le c,
$ and
$$
\bigg\|\frac{1}{n}\sum_{i}\boldsymbol{A}_i(\boldsymbol{\Theta}^\star)\otimes\boldsymbol{x}_i\boldsymbol{x}_i^\top - \boldsymbol{H}^\star\bigg\|_{\max} \leq c\sqrt{\frac{\log(Kd)}{n}},
$$
then $\mathbb{P}(\mathcal{E}_2) \leq c_1 (Kd)^{-c_2}$.

Since the third-order derivatives of the softmax function are bounded, the following deterministic bound holds for any $\boldsymbol{\Theta}$:
$$
\big|[\boldsymbol A_i(\boldsymbol\Theta)]_{kk'}-[\boldsymbol A_i(\boldsymbol\Theta^\star)]_{kk'}\big|\le c_3\sum_{k''}\big|\langle\boldsymbol\theta_{k''}-\boldsymbol\theta^\star_{k''},\boldsymbol x_i\rangle\big|\le c_3\|\boldsymbol x_i\|_\infty\,\|\boldsymbol\Theta-\boldsymbol\Theta^\star\|_{1,1}.
$$
We further derive
$$
\bigg\|\nabla^2\widehat{\mathcal L}_1(\boldsymbol\Theta^\star+t\boldsymbol\Delta)-\nabla^2\widehat{\mathcal L}_1(\boldsymbol\Theta^\star)\bigg\|_{\max}\lesssim\|\boldsymbol\Delta\|_{1,1}.
$$
Since $\check{\boldsymbol\Theta}$ solves the constrained subproblem $\min_{\boldsymbol\Theta:\boldsymbol\Theta_{S^c}=\mathbf{0}}\big\{\widehat{\mathcal L}(\boldsymbol\Theta)+\lambda\|\boldsymbol\Theta\|_{1,1}\big\}$, while $\boldsymbol\Theta^\star$ is also supported on $S$ and hence a feasible solution to this constrained problem. Repeating the argument from the proof of Theorem~\ref{theorem1} yields
$$
\|\check{\boldsymbol\Delta}\|_F\le \frac{9K^2}{4\kappa}\sqrt{Ks}\,\lambda =\frac{9}{4\kappa}K^{5/2}\sqrt{s}\,\lambda. 
$$
Let $\mathcal E_3$ denote the above event, then $\mathbb{P}(\mathcal E_3^c)\le c_3(Kd)^{-1}$.
The vector $\check{\boldsymbol{\Delta}}$ has at most $|S| \leq (K-1)s \leq Ks$ non-zero entries, then
$
\|\check{\boldsymbol{\Delta}}\|_{1,1} = \|\check{\boldsymbol{\Delta}}_S\|_{1,1} \leq \sqrt{|S|}\,\|\check{\boldsymbol{\Delta}}\|_F \leq \sqrt{Ks}\,\|\check{\boldsymbol{\Delta}}\|_F.
$
We obtain 
$$
\|\boldsymbol E\|_{\max}\ \lesssim\ \sqrt{\tfrac{\log(Kd)}{n}}+K^3 s\,\frac{\log(Kd)}{\sqrt n}\ \lesssim\ K^3 s\,\frac{\log(Kd)}{\sqrt n}.
$$
Therefore, when $n\gtrsim\gamma^{-2}\operatorname{poly}(K)s^2\log(Kd)$, we have
$$
D_\infty\,|S|\,\|\boldsymbol E\|_{\max}\le\frac{\gamma}{16}.
$$
We then decompose $\bar{\boldsymbol H}_{SS}$ as
$$
\bar{\boldsymbol H}_{SS}=\boldsymbol H^\star_{SS}+\boldsymbol E_{SS}=\boldsymbol H^\star_{SS}\big(\boldsymbol I+(\boldsymbol H^\star_{SS})^{-1}\boldsymbol E_{SS}\big).
$$
Let $\boldsymbol\Phi:=(\boldsymbol H^\star_{SS})^{-1}\boldsymbol E_{SS}$. Taking the operator $\ell_\infty$ norm yields
$$
|||\boldsymbol\Phi|||_\infty\le|||(\boldsymbol H^\star_{SS})^{-1}|||_\infty\cdot|||\boldsymbol E_{SS}|||_\infty\le D_\infty\cdot|S|\,\|\boldsymbol E\|_{\max}\le\frac{\gamma}{16}\le\frac12.
$$
Hence the Neumann series converges absolutely under the operator $\ell_\infty$ norm, which implies
$$
|||(\boldsymbol I+\boldsymbol\Phi)^{-1}|||_\infty\le\frac{1}{1-|||\boldsymbol\Phi|||_\infty}\le 2.
$$
As a result, $\bar{\boldsymbol H}_{SS}=\boldsymbol H^\star_{SS}(\boldsymbol I+\boldsymbol\Phi)$ is invertible, and its inverse satisfies
\begin{equation}\label{eq:HSSinv}
\begin{aligned}
|||\bar{\boldsymbol H}_{SS}^{-1}|||_\infty
&=|||(\boldsymbol I+\boldsymbol\Phi)^{-1}(\boldsymbol H^\star_{SS})^{-1}|||_\infty\\
&\le|||(\boldsymbol I+\boldsymbol\Phi)^{-1}|||_\infty\cdot|||(\boldsymbol H^\star_{SS})^{-1}|||_\infty\le 2\,D_\infty.
\end{aligned}
\end{equation}

\noindent\textbf{Step 3: Incoherence Condition for the Empirical Hessian.}
Define $\boldsymbol M^\star:=\boldsymbol H^\star_{S^cS}(\boldsymbol H^\star_{SS})^{-1}$, where Assumption~\ref{a2cond1} states $|||\boldsymbol M^\star|||_\infty\le 1-\gamma$. We employ the identity
$$
\bar{\boldsymbol H}_{SS}^{-1}-(\boldsymbol H^\star_{SS})^{-1}=-(\boldsymbol H^\star_{SS})^{-1}\boldsymbol E_{SS}\bar{\boldsymbol H}_{SS}^{-1},
$$
which yields
$$
\begin{aligned}
\bar{\boldsymbol{H}}_{S^c S}\bar{\boldsymbol{H}}_{SS}^{-1} - \boldsymbol{M}^\star &= \boldsymbol{E}_{S^c S}\bar{\boldsymbol{H}}_{SS}^{-1} + \boldsymbol{H}_{S^c S}^\star\big(\bar{\boldsymbol{H}}_{SS}^{-1} - (\boldsymbol{H}_{SS}^\star)^{-1}\big) \\
&= \boldsymbol{E}_{S^c S}\bar{\boldsymbol{H}}_{SS}^{-1} + \boldsymbol{H}_{S^c S}^\star\Big(-(\boldsymbol{H}_{SS}^\star)^{-1}\boldsymbol{E}_{SS}\bar{\boldsymbol{H}}_{SS}^{-1}\Big) \\
&= \boldsymbol{E}_{S^c S}\bar{\boldsymbol{H}}_{SS}^{-1} - \boldsymbol{H}_{S^c S}^\star(\boldsymbol{H}_{SS}^\star)^{-1}\boldsymbol{E}_{SS}\bar{\boldsymbol{H}}_{SS}^{-1} \\
&= \big(\boldsymbol{E}_{S^c S} - \boldsymbol{M}^\star\boldsymbol{E}_{SS}\big)\bar{\boldsymbol{H}}_{SS}^{-1}.
\end{aligned}
$$
Taking the operator $\ell_\infty$ norm and applying the triangle inequality and submultiplicativity of matrix norms, we obtain
$$
\begin{aligned}
||| \bar{\boldsymbol{H}}_{S^c S}\bar{\boldsymbol{H}}_{SS}^{-1} - \boldsymbol{M}^\star |||_\infty &= ||| \big(\boldsymbol{E}_{S^c S} - \boldsymbol{M}^\star\boldsymbol{E}_{SS}\big)\bar{\boldsymbol{H}}_{SS}^{-1} |||_\infty \\
&\le \Big(||| \boldsymbol{E}_{S^c S} |||_\infty + ||| \boldsymbol{M}^\star |||_\infty \cdot ||| \boldsymbol{E}_{SS} |||_\infty\Big) \cdot ||| \bar{\boldsymbol{H}}_{SS}^{-1} |||_\infty \\
&\le \Big(|S|\|\boldsymbol{E}\|_{\max} + 1\cdot |S|\|\boldsymbol{E}\|_{\max}\Big) \cdot 2D_\infty \\
&= 4 D_\infty |S| \|\boldsymbol{E}\|_{\max} \\
&\le \frac{\gamma}{4}.
\end{aligned}
$$
Combining the triangle inequality with Assumption~\ref{a2cond1}, we have
\begin{equation}\label{eq:incoherenceEmp}
|||\bar{\boldsymbol H}_{S^cS}\bar{\boldsymbol H}_{SS}^{-1}|||_\infty\le|||\boldsymbol M^\star|||_\infty+\frac{\gamma}{4}\le(1-\gamma)+\frac{\gamma}{4}=1-\frac{3\gamma}{4}.
\end{equation}

\noindent\textbf{Step 4: Strict Dual Feasibility.}
Take the operator $\ell_\infty$ norm of equation~\eqref{eq:ZSc}. Combining~\eqref{eq:incoherenceEmp} and the bound $\|\check{\boldsymbol Z}_S\|_\infty\le 1$, we have
$$
\|\bar{\boldsymbol H}_{S^cS}\bar{\boldsymbol H}_{SS}^{-1}\check{\boldsymbol Z}_S\|_\infty\le|||\bar{\boldsymbol H}_{S^cS}\bar{\boldsymbol H}_{SS}^{-1}|||_\infty\cdot\|\check{\boldsymbol Z}_S\|_\infty\le 1-\frac{\gamma}{2}.
$$
Applying the inequality $|||\boldsymbol M\boldsymbol v|||_\infty\le|||\boldsymbol M|||_\infty\|\boldsymbol v\|_\infty$ together with~\eqref{eq:incoherenceEmp}, it follows that
\begin{align*}
&\bigg\|\frac{1}{\lambda}\Big(\bar{\boldsymbol H}_{S^cS}\bar{\boldsymbol H}_{SS}^{-1}\big[\nabla\widehat{\mathcal L}(\boldsymbol\Theta^\star)\big]_S-\big[\nabla\widehat{\mathcal L}(\boldsymbol\Theta^\star)\big]_{S^c}\Big)\bigg\|_\infty\\
&\quad\le\frac{1}{\lambda}\Big(|||\bar{\boldsymbol H}_{S^cS}\bar{\boldsymbol H}_{SS}^{-1}|||_\infty+1\Big)\big\|\nabla\widehat{\mathcal L}(\boldsymbol\Theta^\star)\big\|_{\max}\le\frac{2}{\lambda}\big\|\nabla\widehat{\mathcal L}(\boldsymbol\Theta^\star)\big\|_{\max}.
\end{align*}
Substitute the threshold condition $\lambda\ge\frac{8}{\gamma}\big\|\nabla\widehat{\mathcal L}(\boldsymbol\Theta^\star)\big\|_{\max}$ specified in the theorem and combine the two bounds to obtain
\begin{equation}\label{eq:strictDual}
\|\check{\boldsymbol Z}_{S^c}\|_\infty\le1-\frac{\gamma}{2}+\frac{\gamma}{4}=1-\frac{\gamma}{4}<1.
\end{equation}
This completes the proof of Theorem~\ref{theorem2}.

\noindent\textbf{Step 5: $\theta_{\min}$ Condition for Sign Consistency.}
From~\eqref{eq:DeltaS},~\eqref{eq:HSSinv} and $\|\check{\boldsymbol Z}_S\|_\infty\le 1$, we have
\begin{equation}\label{eq:DeltaSbound}
\|\check{\boldsymbol\Delta}_S\|_\infty\le|||\bar{\boldsymbol H}_{SS}^{-1}|||_\infty\Big(\big\|\nabla\widehat{\mathcal L}(\boldsymbol\Theta^\star)\big\|_{\max}+\lambda\Big)\le 2D_\infty\Big(\frac{\gamma}{8}\lambda+\lambda\Big)=2D_\infty\Big(1+\frac{\gamma}{8}\Big)\lambda.
\end{equation}
Conclusion (a) of constrained strong convexity guarantees the existence and uniqueness of the solution to the oracle subproblem as well as the well-posedness of equation~\eqref{eq:DeltaS}.
For any index pair $(j,k)\in S$, combining~\eqref{eq:DeltaSbound} and Assumption~\ref{assumption3} yields
$$
|\widehat\theta_{1,jk}-\theta^\star_{jk}|\le \|\check{\boldsymbol\Delta}_S\|_\infty\le 2D_\infty\Big(1+\frac{\gamma}{8}\Big)\lambda<\theta^\star_{\min}.
$$
This implies no missed support entries over $S$ and identical signs $\operatorname{sign}(\widehat\theta_{1,jk})=\operatorname{sign}(\theta^\star_{jk})$. Together with the inclusion $\operatorname{supp}(\widehat{\boldsymbol\Theta}_1)\subseteq S$ established in Step~4, we finally conclude that
$$
\operatorname{supp}(\widehat{\boldsymbol\Theta}_1)=S,\quad \operatorname{sign}(\widehat{\boldsymbol\Theta}_{1,S})=\operatorname{sign}(\boldsymbol\Theta^\star_S).
$$

Note that both Theorem~\ref{theorem2} and Theorem~\ref{theorem2prime} hold on the same high-probability event $\mathcal E_1\cap\mathcal E_2\cap\mathcal E_3$, which has probability at least $1-C(Kd)^{-c}$ for some constants $c,C>0$.

\end{proof}



\bibliographystyle{plainnat}
\bibliography{references}

\end{document}